\documentclass[12pt]{amsart}
\usepackage{geometry}
\usepackage{amsmath,amsthm,amsfonts,amssymb,bm,graphicx, mathrsfs }
\usepackage{verbatim}

\usepackage
{hyperref}
\hypersetup{colorlinks=true,citecolor=blue,linkcolor=blue,urlcolor=blue}

\theoremstyle{plain}
\newtheorem{theorem}{Theorem}
\newtheorem{corollary}[theorem]{Corollary}
\newtheorem{lemma}{Lemma}

\newtheorem{prop}[theorem]{Proposition}
\numberwithin{equation}{section}

\theoremstyle{definition}
\newtheorem{remark}{Remark}

\newcommand{\dd}{{\rm d}}

\usepackage{subcaption}
\usepackage{graphicx}

\makeatletter
\DeclareRobustCommand\widecheck[1]{{\mathpalette\@widecheck{#1}}}
\def\@widecheck#1#2{%
    \setbox\z@\hbox{\m@th$#1#2$}%
    \setbox\tw@\hbox{\m@th$#1%
       \widehat{%
          \vrule\@width\z@\@height\ht\z@
          \vrule\@height\z@\@width\wd\z@}$}%
    \dp\tw@-\ht\z@
    \@tempdima\ht\z@ \advance\@tempdima2\ht\tw@ \divide\@tempdima\thr@@
    \setbox\tw@\hbox{%
       \raise\@tempdima\hbox{\scalebox{1}[-1]{\lower\@tempdima\box
\tw@}}}%
    {\ooalign{\box\tw@ \cr \box\z@}}}
\makeatother

\begin{document}

\author{Christoph Aistleitner}
\address{Graz University of Technology,
  Institute of Analysis and Number Theory,
  Steyrergasse 30/II,
  8010 Graz,
  Austria}
\email{aistleitner@math.tugraz.at}

\author{Valentin Blomer}
\address{University of G\"ottingen, Mathematisches Institut, Bunsenstr.~3-5, D-37073 G\"ottingen, Germany} \email{vblomer@math.uni-goettingen.de}

\author{Maksym Radziwi\l\l}
\address{
  Caltech,
  Department of Mathematics,
  1200 E California Blvd,
  Pasadena, CA, 91125}
\email{maksym.radziwill@gmail.com}

\title{Triple correlation and long gaps in the spectrum of flat tori}
 

\thanks{The first author is supported by the Austrian Science Fund FWF, projects F-5512 and Y-901. The second author was partially supported by a SNF-DFG lead agency grant BL 915/2-2. The third author acknowledges support of a Sloan fellowship. }

\begin{abstract} We evaluate the triple correlation of eigenvalues of the Laplacian on generic flat tori in an averaged sense. As two consequence we show that (a) the limit inferior (resp. limit superior) of the triple correlation is almost surely at most (resp. at least) Poissonian, and that (b) almost all flat tori contain infinitely many gaps in their spectrum that are at least $2.006$ times longer than the average gap.

The significance of the constant $2.006$ lies in the fact that there exist sequences with Poissonian pair correlation and with consecutive spacings bounded uniformly from above by $2$, as we also prove in this paper. Thus our result goes beyond what can be deduced solely from the Poissonian behavior of the pair correlation.
 \end{abstract}


\setcounter{tocdepth}{2}  

\maketitle 

\section{Introduction}
\subsection{The Berry-Tabor conjecture for tori}
A central objective in quantum chaos is the classification of quantum systems according to universal statistical properties.  
We consider the energy spectrum $\lambda_1\leq \lambda_2 \leq\dots$ of a bound Hamiltonian system, such as the geodesic flow on a compact Riemannian manifold, where the energy levels are the eigenvalues of the associated Laplace-Beltrami operator. According to central conjectures \cite{BT, BGS}, one expects a fundamental dichotomy between integrable and chaotic systems. In particular, the local statistics of the sequence, normalized to have mean spacing unity, should imitate those of certain random matrix ensembles. In the former case, the local statistical properties are expected to coincide with those of a sequence of points coming from a Poisson process, at least in generic cases. In particular,  the probability that   a randomly chosen interval $[X, X+r]$ of fixed length $r$ contains exactly  $k$ of the numbers $\Lambda_j$ is conjectured to be $(k!)^{-1} r^k e^{-r}.$ From a rigorous mathematical standpoint these conjectures are completely out of reach, and we content ourselves here with the simple, yet fundamental example of the following class of  integrable Hamiltonian systems. \\

Let $\Lambda \subseteq \Bbb{R}^2$ be a lattice of rank $2$. Then $\Bbb{R}^2/\Lambda$ endowed with the Euclidean metric is a flat torus and $\Delta = - \frac{\partial^2}{\partial x^2} - \frac{\partial^2}{\partial y^2}$ is the Laplace-Beltrami operator on $\Bbb{R}^2 / \Lambda$. The eigenvalues of the Laplacian $\Delta$ on $\Bbb{R}^2 / \Lambda$ are given by $4 \pi^2 \| {\bm \omega} \|_{2}^{2}$ with ${\bm \omega}$ belonging to the dual lattice $\Lambda^{\star}$. Thus the eigenvalues are integral values of quadratic forms, 
\begin{equation} \label{eq:torii}
q_{\bm \alpha}(m,n) = \alpha_1 m^2 + \alpha_2 m n + \alpha_3 n^2, \quad \bm \alpha = (\alpha_1, \alpha_2, \alpha_3), 
\end{equation}
where $\bm \alpha $ is  always constrained by $4 \alpha_1 \alpha_3 > \alpha_2^2$ since $q_{\bm \alpha}$ is necessarily positive-definite, and ${\bm \alpha}$ is determined by the lattice $\Lambda$ up to ${\rm GL}_2(\Bbb{Z})$-invariance. A standard fundamental domain (up to boundary) is $$\mathcal{D} := \{(\alpha_1, \alpha_2, \alpha_3) \in \Bbb{R}^3 \mid 0 \leq \alpha_2 \leq \alpha_1 \leq \alpha_3\},$$ and  for convenience we will restrict ${\bm \alpha}$ to this set. 
Each quadratic form $q_{\bm \alpha}$ has the automorphism $(m, n) \mapsto (-m, -n)$, so each positive eigenvalue occurs with multiplicity at least 2. It is therefore natural to desymmetrize the spectrum and consider only the values $q_{\bm \alpha}(m, n)$ with $m > 0$ or $m=0$ and $n \geq 0$.  


Given $\bm \alpha = (\alpha_1, \alpha_2, \alpha_3)$ we denote by $0 < \Lambda_1 \leq \Lambda_2 \leq \ldots $ the ordered set of values 
\begin{equation}\label{def-Lambda}
\frac{\pi}{\sqrt{4 \alpha_1 \alpha_3 - \alpha_2^2}} \cdot q_{\bm \alpha}(m,n) 
\end{equation}
where $m > 0$ or $m = 0$ and $n \geq 0$.
Asymptotically the average spacing between the $\Lambda_i$ is one, thus we think of $\Lambda_i$ as the properly re-scaled multi-set of  eigenvalues of the  Laplacian on $\Bbb{R}^2 / \Lambda$ (desymmetrized after removing the obvious double multiplicity   of eigenvalues). 

\subsection{Triple correlation}
In this setting, a consequence of the Berry-Tabor conjecture is that for ``generic'' tori the distribution of the gaps $\Lambda_{i + 1} - \Lambda_{i}$ should coincide with the distribution of nearest neighbor spacings of the Poisson process.\footnote{As an aside, the parallelogram rule implies that there are linear relations in the sequence of $\Lambda_i$, but they do no seem to have an influence on the fine-scale statistical properties.} Precisely, the number of $i \leq N$ for which $\Lambda_{i + 1} - \Lambda_{i} \in [0, r]$ should be asymptotically $N \int_{0}^{r} e^{-u} {\rm d}u$.
The word ``generic'' implicitely assumes a measure, and in this case 
the natural measure on $\mathcal{D}$ 
is the hyperbolic measure
$$
{\rm d} \mu_{\text{hyp}}(\alpha_1, \alpha_2, \alpha_3) = \frac{{\rm d} \alpha_1\, {\rm d} \alpha_2 \, {\rm d} \alpha_3}{(4 \alpha_1  \alpha_3 - \alpha_2^2)^{3/2}}. 
$$



While the Berry-Tabor conjecture for generic flat tori remains wide open, some results have been obtained for \textit{n-correlations}
\begin{displaymath}
\begin{split}
T_{n}(\bm \alpha; I_2, \ldots, I_{n }; N) := \frac{1}{N} \Big| \Big \{(i_1, \ldots, i_n) \in \mathcal{J}_n(N)  \mid  \Lambda_{i_{j}}(\bm \alpha) - \Lambda_{i_{1}}(\bm \alpha) \in I_{j } \text{ for } 2 \leq j \leq n   \Big \}\Big|
\end{split}
\end{displaymath}
where $\mathcal{J}_n(N) \subset [1,N]^{n}$ is the subset of tuples all of whose entries are pairwise distinct and $I_2, \ldots, I_{n}$ are intervals in $\mathbb{R}$. When $I_2 = \ldots = I_{n } = I$ we  write $T_{n}(\bm \alpha; I; N)$. We expect that the $n$-correlations coincide for all $n \geq 1$ with those of a sequence generated by a Poisson process, i.e.\ for almost all $\bm \alpha$ with respect to the measure $\mu_{\text{hyp}}$ and for any fixed intervals $I_2, \ldots, I_n$ we expect to have 
\begin{equation} \label{eq:pair} 
T_{n}(\bm \alpha; I_2, \ldots, I_n; N) \rightarrow \mu(I_2)\cdot \ldots\cdot  \mu(I_{n}) \quad  (N \rightarrow \infty),
\end{equation}
where $\mu$ denotes Lebesgue measure. 

Sarnak \cite{Sa} established \eqref{eq:pair} for $n = 2$ for almost all $\bm \alpha$ with respect to the hyperbolic measure $ \mu_{\text{hyp}}$ (equivalently for almost all $\bm \alpha$ with Lebesgue measure on $\alpha_1, \alpha_2, \alpha_3$). This was generalized in \cite{Va} to four-dimensional tori. Eskin, Margulis and Mozes \cite{EMM} gave explicit Diophantine conditions on ${\bm \alpha}$ under which \eqref{eq:pair} holds when $n = 2$.

In contrast nothing is known for the case of correlations of higher order, that is, for $n > 2$. While the pair correlation is associated to an orthogonal group of signature $(2,2)$, no useful group structure seems to be available in the case of higher correlations. In particular, the methods of \cite{EMM} are not available. Our main result sheds some light on the case $n > 2$ by establishing \eqref{eq:pair} for $n = 3$ in an averaged sense.

\begin{theorem}\label{thm:main} Let $R \subseteq \mathcal{D}$ be a  three-dimensional rectangle of finite hyperbolic measure $\mu_{\text{{\rm hyp}}}(R)$, and let $J \subset \mathbb{R}$ be a finite interval.  Then
\begin{equation}\label{triple}
 \lim_{N \rightarrow \infty} \int_{R}  T_3({\bm \alpha}; I_1, I_2; N) \,\dd \mu_{\text{{\rm hyp}}}(\bm \alpha) = \mu(I_1) \mu(I_2) \mu_{\text{{\rm hyp}}}(R) ,
\end{equation}
uniformly for all intervals $I_1, I_2 \subseteq J$. 
\end{theorem}
A consequence of Theorem \ref{thm:main} is a sharp upper bound for $\liminf T_{3}(\bm \alpha; I; N)$ for almost all $\bm \alpha \in \mathcal{D}$ with respect to $\mu_{\text{hyp}}$. 

\begin{corollary}\label{cor:liminf} Let $I \subseteq \mathbb{R}$ be a finite interval. For almost all ${\bm \alpha} \in \mathcal{D}$ with respect to the measure $\mu_{\text{{\rm hyp}}}(\bm \alpha)$ we have
\begin{equation*}
  \liminf_{N \rightarrow \infty} T_3({\bm \alpha}; I; N) \leq \mu(I)^2 \quad \text{and} \quad \limsup_{N \rightarrow \infty} T_{3}({\bm \alpha}; I; N) \geq \mu(I)^2.
\end{equation*}
\end{corollary}

Although the statistical investigation of the spectrum of a Riemann surface is a problem of mathematical physics, the main input of the proof of Theorem \ref{thm:main} is of Diophantine nature and consists in an involved lattice point count in a certain 6-dimensional region, see Section \ref{sec:Dioph}. It uses to some extent ideas of \cite{Sa}, but   since Sarnak studies an $L^2$-mean of the pair correlation, he has 8 variables at his disposal, while we have only 6. Consequently, the Diophantine analysis becomes more delicate and depends, among other things, on non-trivial bounds for Kloosterman sums. In fact our result can be viewed as a ``thin subset'' version of Sarnak's result. These remarks are explained in more detail in Section \ref{sec:thmain}.

\subsection{Long gaps in Poissonian sequences} Theorem \ref{thm:main} was originally motivated by the problem of establishing the existence of large gaps between consecutive $\Lambda_{i}$. This problem stands in sharp contrast with the problem of establishing small gaps between consecutive $\Lambda_{i}$.  In the case of rectangular tori (i.e.\ $\alpha_2 = 0$) it is possible to obtain nearly optimal results for almost all $\bm {\alpha}$, that is, gaps of size $N^{-1 + \varepsilon}$ among the first $N$ eigenvalues (birthday paradox), see \cite{BBRR}. As long gaps in a Poisson process typically grow at a logarithmic scale, their investigation is a more delicate endeavour. In addition, long gaps are harder to detect from a technical point of view. While a difference of any two eigenvalues is in particular an upper bound for the smallest gap, no such approximation  is available for a lower bound of long gaps.


It is therefore not surprising that  comparatively little is known about large gaps between consecutive $\Lambda_i$. Using a variant of Corollary \ref{cor:liminf} (see Lemma \ref{lem4}) we obtain the following partial progress on this question.

\begin{corollary}\label{cor:gaps} Let $G$ be the unique positive solution to the equation
  $$
  35 - 9 G^2  + (12 G - 10) \sqrt{6 G - 5} = 0.
  $$
  Then, for almost all $\bm \alpha$ with respect to the measure $\mu_{\text{{\rm hyp}}}$ we have
  $$
  \limsup_{i \rightarrow \infty} \ (\Lambda_{i + 1}(\bm \alpha) - \Lambda_i(\bm \alpha)) \geq G = 2.00636193892510 \ldots 
  $$
\end{corollary}
\begin{remark}
  In fact we prove a slightly stronger result: Given $\varepsilon > 0$, for almost all $\bm \alpha$ there exists a subsequence $N_1 < N_2 < \ldots$ along which we have $\Lambda_{i + 1} - \Lambda_{i} > G - \varepsilon$ for a positive proportion of $i \leq N_{j}$ as $N_{j} \rightarrow \infty$. 
\end{remark}

As mentioned above the Berry--Tabor conjecture predicts that the gaps $\Lambda_{i + 1}(\bm \alpha) - \Lambda_i(\bm \alpha)$ are exponentially distributed (for generic $\bm \alpha$), which clearly would imply that the maximal gap size is unbounded. However, proving this seems to be far out of reach. As a remark pointing to this direction, note that by inclusion-exclusion it is not hard to see that if we could show \eqref{eq:pair} for all   $n   \leq 2k$ ($k \in \Bbb{N}, k \geq k_0$) and some ${\bm \alpha}$, then the corresponding torus has infinitely many gaps in its spectrum that are of length at least $k/2$ 
In particular, if \eqref{eq:pair} is known for all $n$, one obtains unbounded gaps.

The particular significance of Corollary \ref{cor:gaps} lies in the fact that the spacings that we obtain are strictly greater than $2$. As we show in the following theorem, there exist sequences with gaps uniformly bounded by $2$ and whose pair correlation is Poissonian. In other words, Theorem \ref{thm:main} shows that the gap structure of almost all flat tori has properties that cannot be derived from the pair correlation result alone. 

It is an interesting problem in its own right to determine the minimal possible value of the largest gap of a sequence whose pair correlation is assumed to be Poissonian.  To fix notation, let us denote by $\mathscr{S}$ the set of all increasing sequences $(\gamma_i)$ of non-negative real numbers with mean-spacing one and Poissonian pair correlation, i.e.\
\begin{equation}\label{sequence}
  \frac{1}{N} \gamma_N \rightarrow 1\quad \text{and} \quad \frac{1}{N} \sum_{\substack{\gamma_{i_1} - \gamma_{i_2} \in I \\ i_1, i_2 \leq N,~i_1 \neq i_2}} 1 \rightarrow \mu(I) \qquad \text{($N \rightarrow \infty$)}
\end{equation}
for all finite intervals $I \subset \mathbb{R}$.

\begin{theorem}\label{thm4} Let  $$\nu := \inf_{(\gamma_i) \in \mathscr{S}} \limsup_i (\gamma_{i+1} - \gamma_i).$$ Then $\nu \in [3/2, 2].$
\end{theorem}



The lower bound $\nu \geq 3/2$ shows that every sequence with Poissonian pair correlation contains infinitely many gaps that are 1.5 times larger than the average gap. By \cite{EMM}, this applies for instance to the ordered sequence of values $x^2 + \sqrt{2} y^2$, $x, y \in \Bbb{N}$.  The lower bound is obtained using a reasoning similar to the one which leads from Corollary \ref{cor:liminf} to Corollary \ref{cor:gaps}. 

The problem of determining $\nu$ precisely is interesting in its own right since it touches upon the questions of how one can construct a sequence having a preassigned pair correlation distribution, and for which pair correlation distributions such a construction can exist at all and for which it cannot. We will discuss this problem a bit further in the closing remarks of the introduction (Subsection \ref{subs:close}).

Finally we remark that there are interesting sequences whose pair correlation is Poissonian, but whose level spacing is not, for instance the (suitably normalized) sequence of fractional parts of $\sqrt{n}$, cf.\ \cite{EM, EMa}, which has in fact gaps of unbounded length.

\subsection{Outline of the proof of Theorem \ref{thm:main}} \label{sec:thmain}

 Roughly speaking, establishing Theorem \ref{thm:main} amounts to evaluating asymptotically 
 $$
\frac{1}{N} \int_{(\alpha, \beta) \in \mathcal{R}} \int_{N}^{2N} \Big ( \sum_{\substack{|u - \alpha m^2 - m n - \beta n^2| \leq \eta \\ M \leq m, n \leq 2M}} 1 \Big )^3 {\rm d} u \, {\rm d} \alpha \, {\rm d} \beta
 $$
as $N \rightarrow \infty$, 
 with $\eta > 0$ fixed, $M^2 = N$ and $\mathcal{R}$ a fixed rectangle.   Note that for $\alpha_2 \not= 0$ we have $\Lambda_i ((\alpha_1, \alpha_2, \alpha_3)) = \Lambda_i  (  (  \alpha_1/\alpha_2, 1 ,  \alpha_3/\alpha_2   ) )$, so without loss of generality we can assume $\alpha_2 = 1$.  
 Expanding the third power gives rise to three conditions $|u - \alpha m_i^2 - m_i n_i - \beta n_i^2| \leq \eta$ with $i = 1,2,3$ and $m_i, n_i \in [M, 2M]$ say. We detect each of them using a Fourier integral. Subsequently we execute the integrations in $u, \alpha, \beta$ and make the following change of variables
 \begin{align*}
   & a_1 = n_1 - n_2, \ a_2 = n_1 + n_2, \ a_3 = m_1 - m_3, \ a_4 = m_1 + m_3 ,\\
   & b_1 = m_1 - m_2, \ b_2 = m_1 + m_2, \ b_3 = n_1 - n_3, \ b_4 = n_1 + n_3 .
 \end{align*}
 Here the new variables are constrained by $a_1 + a_2 = b_3 + b_4$ and $b_1 + b_2 = a_3 + a_4$. 
 After this change of variables and repeated integration by parts in the resulting Fourier integral, we are lead to the problem of understanding asymptotically an expression that roughly looks like
 \begin{equation} \label{eq:todo}
 \sum_{\substack{M \leq a_i, b_i \leq 2M \\ |\Delta_1|, |\Delta_2| \ll \Delta + M^2 \\ a_1 + a_2 = b_3 + b_4 \\ b_1 + b_2 = a_3 + a_4}}  \frac{1}{\Delta + M^2}  
 \end{equation}
 where
 $$
 \Delta = a_1 a_2 a_3 a_4 - b_1 b_2 b_3 b_4
 $$
 and
 \begin{align*}
   & \Delta_1 = a_1 a_2 b_3 a_4 + a_1 a_2 a_3 b_4 - a_1 b_2 b_3 b_4 - b_1 a_2 b_3 b_4 ,\\
   & \Delta_2 = a_1 b_2 a_3 a_4 + b_1 a_2 a_3 a_4 - b_1 b_2 a_3 b_4 - b_1 b_2 b_3 a_4. 
 \end{align*}
 The expected main term for \eqref{eq:todo} is of size $M^2$.
 It is important to note that we are oversimplifying the situation here by assuming that $a_i, b_i$ are of size $M$. In practice this is not necessarily the case, but we assume it here for the sake of simplicity.

 In \cite{Sa}, Sarnak obtains the same expression as \eqref{eq:todo} but without the constraints $a_1 + a_2 = b_3 + b_4$ and $b_1 + b_2 = a_3 + a_4$. As a result his expression is asymptotically of size $M^4$ and the counting problem is easier since there are fewer constraints. This explains our earlier remark on our work being a ``thin subset'' version of Sarnak's work. 
 
 When $|\Delta|  > M^{4 - \delta}$ we can estimate \eqref{eq:todo} asymptotically by applying the Lipschitz principle (that is, the Euler-Maclaurin approximation), getting a main term of size $M^2$. It therefore remains to show that the contribution to \eqref{eq:todo} of the terms with $|\Delta| \leq M^{4 - \delta}$ is negligible, that is $\ll M^{2 - \varepsilon}$ for some $\varepsilon > 0$. 

 In the case of $|\Delta| \in [D, 2D]$ with $D < M^{4 - \delta}$ we use a substitution trick of Sarnak to see that the condition $|\Delta_1| \ll D + M^2$ implies that $|k l| \ll D + M^2$ where
 \begin{equation} \label{eq:deff}
 k := a_2 a_4 - b_2 b_4 , \ l = b_2 b_3 - a_2 a_3.
 \end{equation}
 In Sarnak's case this alone is enough to conclude via an elementary argument. In our case his elementary argument barely fails to be sufficient, and we need to non-trivially exploit the new condition that $a_1 + a_2 = b_3 + b_4$ and $a_3 + a_4 = b_1 + b_2$.

 We proceed by
 \begin{enumerate}
 \item\label{1} solving for $a_4$ in terms of $a_2, a_3, b_2, k, l$ using the conditions $a_1 + a_2 = b_3 + b_4$ and $a_3 + a_4 = b_1 + b_2$, Definition \eqref{eq:deff}, and the condition $|\Delta| \in [D, 2D]$. We refer to this as the ``geometric restriction'' on $a_4$.
 \item expressing $b_3$ in terms of the variables $a_2, a_3, b_2, l$ by requiring that $a_3 \equiv - \overline{a_2} l \pmod{b_2}$ and
   $$
 (b_3 = )  \frac{a_3 a_2 + l}{b_2} \approx M.
   $$
 \item expressing $b_4$ in terms of the variables $a_2, a_4, b_2, k$ by requiring that $a_4 \equiv \overline{a_2} k \pmod{b_2}$ and
   $$
 (b_4 = )  \frac{a_2 a_4 - k}{b_2} \approx M. 
   $$
 \end{enumerate}
 This leads to an upper bound for $\eqref{eq:todo}$ that is roughly of the form
 \begin{equation} \label{outline}
 \frac{1}{D + M^2} \sum_{a_2, b_2 \approx M} \sum_{k,l \approx M^2} \sum_{\substack{a_3 \equiv - \overline{a_2} l  \, (\text{mod } b_2)\\ (a_3 a_2 + l) / b_2 \approx M}} \sum_{\substack{a_4 \equiv \overline{a_2} k \, (\text{mod } b_2) \\ (a_2 a_4 - k) / b_2 \approx M}} \Phi_{a_2, a_3, b_2, k, l}(a_4),
   \end{equation}
   where $\Phi$ is a smooth function capturing the ``geometric restriction'' on $a_4$ in terms of the variables $a_2, a_3, b_2, k, l$ that we mentioned in $(1)$ of the preceding list. The function $\Phi$ needs only to capture a barely non-trivial piece of this ``geometric restriction''. It is not necessary for $\Phi$ to be as precise so as to allow to exactly reconstruct $a_4$ in terms of the other variables $a_2, a_3, b_2, k, l$.  

   We now apply Poisson summation on the variables $a_2, a_4$. The diagonal term gives  a negligible contribution thanks to the presence of $\Phi$ (if $\Phi$ were replaced by $1$ we would have obtained a diagonal of size matching the main term $M^2$ and this would not allow us to win). The off-diagonal terms give rise to sums of Kloosterman sums and we win by applying the Weil bound. It is, however, important to design $\Phi$ very carefully as to  capture only a barely non-trivial piece of our condition \eqref{1}, otherwise we encounter again a counting problem with highly cuspidal regions which is hard to analyze. 

  We close by recalling that this discussion assumed that the variables $a_i, b_i$ are in generic position, that is, all of size about $M$. However ranges in which $a_i$ are smaller require different treatement. For instance in some ranges we continue after \eqref{outline} by also applying Poisson summation in $k$ (and even $l$). Moreover in the non-generic ranges the diagonal gives an acceptable contribution and we do not need to introduce the function $\Phi$. This is welcome, since the behavior of $\Phi$ is complicated so we are happy to avoid it whenever possible. 

  \subsection{Outline of the proof of Theorem \ref{thm4}}

  The lower bound in Theorem \ref{thm4} comes from noticing that if the gaps $\gamma_{i + 1} - \gamma_{i}$ are uniformly bounded by, say $3/2 - \varepsilon$, then the distribution function of the gaps $\gamma_{i + 1} - \gamma_{i}$ can grow at most linearly between $1/2  - \varepsilon$ and $3/2  - \varepsilon$ because the pair correlation function is Poissonian. However, integrating by parts we see that there can be no such distribution function unless $\varepsilon = 0$, since the average spacing of the sequence is 1 by assumption.

  The construction of the upper bound is more involved and proceeds as follows. In each interval $[i, i + 1]$ we first select a point $\gamma_{i}$ uniformly at random. The resulting sequence is such that almost surely
  \begin{equation} \label{almost-pair}
  \frac{1}{N} \sum_{\substack{\gamma_{i} - \gamma_{j} \leq b \\ i > j}} 1 \rightarrow \begin{cases}
    \frac{b^2}{2} & \text{ if } b \in [0,1], \\
    b - \frac{1}{2} & \text{ if } b \geq 1.
  \end{cases}
\end{equation}
and the mean spacing of $\gamma_{i}$ is $1$ by construction.

We then make a deterministic correction to the above random construction by introducing a small amount of clusters of points. Since there are few clusters and they are far apart they will not interact with our random construction. Specifically for each $m < \sqrt{N}$ we insert into each interval $[m^2 , m^2 + 1]$ a sequence of $\lceil  \sqrt{2 m} \rceil$ equally spaced points. For $b \geq 1$ each such interval contribute an additional $\frac{1}{2} \cdot \lceil \sqrt{2 m} \rceil^2$ to \eqref{almost-pair}, and therefore a total of $\frac{1}{2}$ when summed over all $m < \sqrt{N}$. On the other hand, for $b \leq 1$ each interval $[m^2 , m^2 + 1]$ contributes an additional
  $$
  \frac{1}{N} \cdot \Big ( \lceil \sqrt{2 m} \rceil (1 - b) \cdot \lceil \sqrt{2m} \rceil b + \sum_{k < b \lceil \sqrt{2 m} \rceil} k \Big ) \sim \frac{2m}{N} \cdot \Big ( b - \frac{b^2}{2} \Big )
  $$
  and summing over $m < \sqrt{N}$ this adds exactly the missing $b - \frac{b^2}{2}$ to \eqref{almost-pair} when $b \leq 1$. 

  \subsection{Closing remarks} \label{subs:close}

  We have not been able to construct a sequence with bounded gaps whose pair correlation and triple correlation are both Poissonian; it would be interesting to know if such a sequence exists, and which bound for the maximal gap size can be achieved. More generally it is an interesting question to ask which distribution functions can occur for the distribution function of the spacings $\gamma_{i + 1} - \gamma_{i}$ given that the pair correlation function of $\gamma_{i}$ is Poissonian. Concerning the general problem of generating a sequence having preassigned pair correlation behavior, there exist several result in this direction in the case when the sequence is given by a so-called random point processes; see for example \cite{kls1,kls2} and the references given there. However, these results only cover the purely random case. As far as we know, the question which asymptotic distributions of the pair correlation can be realized by a deterministic sequence $(\gamma_i)$ has never been studied.  We hope that our investigations initiate further research in this direction.\\

On a somewhat unrelated theme, we note that if we attempt to compute the $4$-correlations, then the first problem we face is that we have to estimate the number of solutions in $M \leq x_1, x_2, x_3, x_4, y_1, y_2, y_3, y_4 \leq 2M$ to
   \begin{equation} \label{eq:var}
   -M^4 \leq 
   \det
   \begin{vmatrix}
     x_1^2 - x_2^2 & x_1^2 - x_3^2 & x_1^2 - x_4^2 \\
     y_1^2 - y_2^2 & y_1^2 - y_3^2 & y_1^2 - y_4^2 \\
     x_1 y_1 - x_2 y_2 & x_1 y_1 - x_3 y_3 & x_1 y_1 - x_4 y_4 
   \end{vmatrix}
   \leq M^{4},
   \end{equation}
    asymptotically as $M \rightarrow \infty$. At present we do not know how to accomplish this, since this corresponds to a very thin region. What complicates matters is that the variety in \eqref{eq:var} is highly singular. If it were smooth then the methods of \cite{JingJing} might offer some hope. 

   Any attempt to estimate the $3$-correlations in an $L^2$ sense runs into a similar problem. Precisely we would have to asymptotically estimate as $M \rightarrow \infty$ the number of solutions in $M \leq m_1, m_2, m_3, m_1', m_2', m_3', n_1, n_2, n_3, n_1', n_2', n_3' \leq 2M$ to
   $$
   -D \leq \det \begin{vmatrix}
     m_1^2 - m_2^2 & m_1^2 - m_3^2 & m_2'^2 - m_1'^2 & m_3'^2 - m_1'^2 \\
     n_1^2 - n_2^2 & n_1^2 - n_3^2 & n_3'^2 - n_1'^2 & n_3'^2 - n_1'^2 \\
     m_1 n_1 - m_2 n_2 & m_1 n_1 - m_3 n_3 & 0 & 0 \\
     0 & 0 & m_2' n_2' - m_1' n_1' & m_3' n_3' - m_1' n_1' 
     \end{vmatrix} \leq D
     $$
     uniformly in $M^{6} \leq D \leq M^{8}$ and the case $D = M^{6}$ appears to be of comparable difficulty to \eqref{eq:var}.


     \subsection{Notational Conventions}

For the rest of the paper, we will apply the following {notational conventions}. Boldface letters like $\textbf{a}$ or $\textbf{A}$ denote  vectors with components $(a_1, \ldots, a_n)$ resp.\ $(A_1, \ldots, A_n)$ whose dimension $n$ will always be clear from the context.  Given a vector $\textbf{v} = (v_1, \ldots, v_n)$ and a real $M \geq 0$, the notation $\textbf{v} \ll M$ means that there exists an absolute constant $C > 0$ such that $|v_i| \leq C M$ for all $i = 1,\ldots, n$. In Section \ref{sec:Dioph} we will make use of the notation $X \preccurlyeq Y$ to mean $X \ll_{\varepsilon} Y M^{\varepsilon}$ for $\varepsilon > 0$ where the meaning of $\varepsilon$ can change from line to line. 
 
 Given a Schwartz function $h$ we define its Fourier transform as
 $$
 \hat{h}(\xi_1, \ldots, \xi_n) := \int_{\mathbb{R}^{n}} h(x_1, \ldots, x_n) e(-x_1 \xi_1 - \ldots - x_n \xi_n) \dd \xi_1 \ldots \dd \xi_{n},
 $$
 where $e(x) := e^{2\pi i x}$. For two vectors $\textbf{v} = (v_1, \ldots, v_n)$ and $\textbf{w} = (w_1, \ldots, w_n)$ we denote by $\psi(\textbf{v}, \textbf{w})$ the function $\psi(v_1, \ldots, v_n, w_1, \ldots, w_m)$.
 Moreover for  $\textbf{v}=(v_1, \ldots, v_n)$ we write
 $$
 \sum_{\textbf{v}} \psi(\textbf{v}) := \sum_{v_1, \ldots, v_n \in \mathbb{Z}} \psi(v_1, \ldots, v_n) \ , \ \int_{\mathbb{R}^n} \psi(\textbf{x}) \dd \textbf{x} := \int_{\mathbb{R}^n} \psi(x_1, \ldots, x_n) \dd x_1 \ldots \dd x_{n},
 $$
 and we will also sometimes write $ \dd (x_1, \ldots, x_n)$ in place of $\dd x_1 \ldots \dd x_{n}$.  \\

\section{Proof of Theorem \ref{thm:main}: Reduction to Diophantine problems}\label{sec4}

Since $T_{3}(\lambda \bm \alpha; I_1, I_2; N) = T_{3}(\bm \alpha; I_1, I_2; N)$ by \eqref{def-Lambda} for all $\lambda > 0$, we have for $\bm \alpha = (\alpha_1, \alpha_2, \alpha_3)$ and $R = I \times J \times K \subseteq \mathcal{D}$ with $I,J,K$ intervals of $\mathbb{R}_{\geq 0}$ that
 \begin{equation} \label{eq:factor}
 \int_{R} T_{3}(\bm \alpha; I_1, I_2; N) d \mu_{\text{hyp}}(\bm \alpha) = \int_{J} \int_{I / \alpha_2} \int_{K / \alpha_2} T_{3}((\alpha_1, 1, \alpha_3); I_1, I_2; N) \dd_{\text{hyp}} (\alpha_1,1,  \alpha_3)   \frac{\dd \alpha_2}{\alpha_2}   
 \end{equation}
 where
 \begin{equation}\label{measure}
 \dd_{\text{hyp}} (\alpha_1, 1, \alpha_3) = \frac{\dd \alpha_1 \, \dd \alpha_3}{(4 \alpha_1 \alpha_3 - 1)^{3/2}}
 \end{equation}
 and $I/\alpha_2$, $K/\alpha_2$ are the corresponding re-scaled intervals with the interpretation $I/\alpha_2 = \Bbb{R}_{\geq 0}$ for $\alpha_2 = 0$. Since  $\{\bm \alpha  \subseteq \mathcal{D} \mid \alpha_2 = 0\}$ is a null set, to prove Theorem \ref{thm:main} it is enough to obtain an asymptotic formula for
 $$
 \int_{R} T_{3}(\bm \alpha; I_1, I_2; N) \dd_{\text{hyp}} \bm \alpha
 $$
 with $R \subseteq \mathcal{D}_{0} := \{(\alpha_1, \alpha_3) \in (0, \infty)^2 : 1 \leq \alpha_1 \leq \alpha_3\}$
a rectangle, $\bm \alpha = (\alpha_1, 1, \alpha_3)$ and $\dd_{\text{hyp}}(\bm \alpha)$ defined as in \eqref{measure}. This is the purpose of Proposition \ref{prop1} below. As before we write $q_{{\bm \alpha}}(m, n) = \alpha_1 m^2 + mn + \alpha_3 n^2$ for ${\bm \alpha} = (\alpha_1, 1, \alpha_3)$. We also define
 $$D({\bm \alpha}) = \frac{1}{\pi} \sqrt{4 \alpha_1\alpha_3 - 1},$$
 so that $\{q_{{\bm \alpha}}(m, n)/D({\bm \alpha}) \mid m > 0 \text{ or } m = 0, n \geq 0 \}$ is the sequence of $\Lambda_j({\bm \alpha})$.
 
 \begin{prop}\label{prop1} Let $F$, $W_i$ ($i = 1, 2$), and $V$   be fixed smooth functions with compact support in $\mathcal{D}_{0}$, $(-\infty, \infty)$ and $(0, \infty)$ respectively.  For $M \geq 1$ let
 \begin{equation*}
\begin{split}
\mathcal{T}(M) := \frac{1}{8} \int_{\mathcal{D}_0}  F({\bm \alpha})\left.\sum_{\textbf{x}, \textbf{y}} \right.^{\ast} 
& W_1\left(\frac{q_{{\bm \alpha}}(x_1, y_1) - q_{{\bm \alpha}}(x_2, y_2)}{D({\bm \alpha}) }\right)\\
&W_2\left(\frac{q_{{\bm \alpha}}(x_1, y_1) - q_{{\bm \alpha}}(x_3, y_3)}{D({\bm \alpha}) }\right)V\left(\frac{q_{{\bm \alpha}}(x_1, y_1)}{M^2 D({\bm \alpha}) }\right) \, \dd_{\text{{\rm hyp}}}{\bm \alpha},
\end{split}
\end{equation*}
where ${\bm \alpha} = (\alpha_1, 1, \alpha_3)$, $\textbf{x} = (x_1, x_2, x_3) \in \Bbb{Z}^3$, $\textbf{y} = (y_1, y_2, y_3) \in \Bbb{Z}^3$ and $\sum^{\ast}$  indicates the conditions 
\begin{equation}\label{cond}
(x_1,y_1) \not = \pm (x_2, y_2), \quad   (x_1,y_1) \not = \pm (x_3, y_3), \quad (x_2, y_2) \not= \pm (x_3, y_3).
\end{equation}
Then  
$$\mathcal{T}(M) = \text{{\rm vol}}(F) \hat{W}_1(0) \hat{W}_2(0)  \hat{V}(0) M^2 + O\left(M^{2 - \frac{1}{258}}\right),$$
where $$\text{{\rm vol}}(F) = \int_{\mathcal{D}_0} F({\bm \alpha}) \dd_{\text{{\rm hyp}}}({\bm \alpha}), \quad {\bm \alpha} = (\alpha_1, 1, \alpha_3). $$
  \end{prop}

  The factor $\tfrac{1}{8}$ is included to compensate for the fact that we are summing over $\textbf x$ and $\textbf y$ of arbitrary signs, and it  therefore takes care of the desymmetrization of the spectrum mentioned in the  introduction. 
 Theorem \ref{thm:main} is an immediate consequence of Proposition \ref{prop1}: choosing $M^2 = N$ and choosing $F$, $W_i$ and $V$ to be smooth approximations to the characteristic functions on $R$, $I_i$ ($i = 1, 2$) and $[2^{-k-1}, 2^{-k}]$ ($k = 0, 1, 2, \ldots, K = 1/\varepsilon$, say)  respectively, we obtain
\begin{displaymath}
\begin{split}
   \int_R  \frac{1}{N} & \Big|\Big\{(i, j, k) \mid \begin{matrix}   \Lambda_j({\bm \alpha}) - \Lambda_i({\bm \alpha}) \in I_1,  \Lambda_k({\bm \alpha}) - \Lambda_i({\bm \alpha}) \in I_2 \\  |\{i,j, k\}| =3, 1 \leq \Lambda_i({\bm \alpha}) \leq N \end{matrix}\Big\} \Big| \dd_{\text{hyp}}\bm \alpha \sim \mu_{\text{hyp}}(R) 
   \mu(I_1) \mu(I_2).   
\end{split}
\end{displaymath}
as $N \rightarrow \infty$. Plugging this into \eqref{eq:factor} we obtain \eqref{triple}, noting that the condition $N \leq \Lambda_i({\bm \alpha}) \leq 2N$ can be replaced by the condition $N \leq i < 2N$, since the cardinality of the symmetric difference of $\{i : N \leq \Lambda_i({\bm \alpha}) \leq 2 N\}$ and $[N, 2N)$ is $o(N)$ as $N \rightarrow \infty$. The rate of convergence is continuous in the endpoints of $I$, and hence uniform as long as $I$ varies within a fixed finite interval $J$.\\

We start the \textbf{proof of Proposition \ref{prop1}} by Fourier-inverting the weight functions $W_1, W_2$ and $V$. The support of $V, W_1, W_2$ implies automatically $\textbf{x}, \textbf{y} \ll M$. We remember this by inserting a smooth redundant weight function $\psi(\textbf{x}/M, \textbf{y}/M)$, where $\psi$ is a suitable smooth function that is constantly 1 on some sufficiently large fixed box in $\Bbb{R}^6$ and constantly zero outside some slightly larger box. We obtain
\begin{equation}\label{step1}
\begin{split}
\mathcal{T}(M) = & \frac{1}{8} \int_{\mathcal{D}_0} F({\bm \alpha})\underset{\textbf{x}, \textbf{y}} {\left.\sum\right.^{\ast}} \psi\left( \frac{\textbf{x}}{M}, \frac{\textbf{y}}{M}\right)\\
& 
\int_{\Bbb{R}^3} \hat{V}(z) \hat{W}_1(u)\hat{W}_2(v)e\left(u\frac{q_{{\bm \alpha}}(x_1, y_1) - q_{{\bm \alpha}}(x_2, y_2)}{D({\bm \alpha}) }\right) \\
&e\left(v\frac{q_{{\bm \alpha}}(x_1, y_1) - q_{{\bm \alpha}}(x_3, y_3)}{D({\bm \alpha}) }\right) e\left(\frac{z q_{{\bm \alpha}}(x_1, y_1) }{M^2 D({\bm \alpha}) }\right) \,  \dd(u, v, z)\,  \frac{\dd{\bm \alpha}}{\pi^3 D({\bm \alpha})^3}. 
\end{split}
\end{equation}
We change variables $ z\leftarrow z/D({\bm \alpha}) $, $u \leftarrow u/D({\bm \alpha}) $, $v \leftarrow v/D({\bm \alpha}) $ and  integrate over ${\bm \alpha} = (\alpha_1, 1, \alpha_3) \in \mathcal{D}_0 \subseteq \Bbb{R}^2$. This gives
\begin{equation}\label{step2}
\begin{split}
&\mathcal{T}(M) =  \frac{1}{8} \underset{\textbf{x}, \textbf{y} } {\left.\sum\right.^{\ast}}\psi\left( \frac{\textbf{x}}{M}, \frac{\textbf{y}}{M}\right)  
 \int_{\Bbb{R}^3}  e\left(u(x_1y_1 - x_2y_2) + v(x_1y_1 - x_3y_3) + \frac{z}{M^2} x_1y_1\right) \\
 &   \widehat{G}\left(u(x_1^2 - x_2^2) + v(x_1^2 - x_3^2) + \frac{z x_1^2}{M^2}, u(y_1^2 - y_2^2) + v(y_1^2 - y_3^2) + \frac{zy_1^2}{M^2} ; z, u, v\right)    \dd(u, v, z)
   \end{split}
\end{equation}
where $$G(\alpha_1, \alpha_3; z, u, v) = \pi^{-3} F({\bm \alpha})  \hat{V}(z D({\bm \alpha}) )  \hat{W}_1(uD({\bm \alpha}) )\hat{W}_2(vD({\bm \alpha}) ) $$
and the Fourier transform $\widehat{G}$ of $G$ is taken with respect to the first two variables $\alpha_1, \alpha_3$. Note that $G$ is a Schwartz-class function in all variables.  
Therefore $\widehat{G}$ is a Schwartz-class function in all variables as well, i.e.\
\begin{equation}\label{schwartz}
\mathscr{D}\widehat{G}(U, V; z, u, v) \ll_{A, \mathscr{D}} \big((1+|U|)(1 + |V|)(1+|z|)(1 + |u|)(1+|v|)\big)^{-A}
\end{equation}
for all $A > 0$ and any differential operator $\mathscr{D}$.\\

 Next  we make an invertible change of integer variables as follows: each $6$-tuple $(x_1, y_1, x_2, y_2, x_3, y_3)$ satisfying \eqref{cond}   is mapped bijectively to an $8$-tuple $(a_1, \ldots, a_4, b_1, \ldots, b_4)$ satisfying 
\begin{equation}\label{cong}
a_1 \equiv a_2 \, (\text{mod } 2), 
\quad  b_1 \equiv b_2 \, (\text{mod } 2)
\end{equation}
as well as 
\begin{equation}\label{same}
\begin{split} 
 & (a_1, a_2, a_3, a_4) \not= (b_3, b_4, b_1, b_2), \quad (a_1, a_2, a_3, a_4) \not= (b_4, b_3, b_2, b_1),\\
  &(a_1, b_1) \not= (0, 0), \quad (a_2, b_2) \not= (0, 0), \quad (a_3, b_3) \not= (0, 0), \quad (a_4, b_4) \not= (0, 0)\\
 \end{split} 
\end{equation}
and 
\begin{equation}\label{sum}
a_1+a_2 = b_3+b_4, \quad a_3 + a_4 = b_1+b_2
\end{equation}
via 
\begin{equation}\label{step3}
\begin{split}
  & a_1 = y_1 - y_2, \quad a_2 = y_1 + y_2,\quad a_3 = x_1 - x_3,\quad a_4 =x_1 + x_3 ,\\
  &b_1 = x_1 - x_2, \quad b_2 = x_1 + x_2,\quad b_3 = y_1 - y_3,\quad b_4 = y_1 + y_3.\\
\end{split}
\end{equation}
Note that the conditions  $a_3 \equiv a_4$  (mod  $2$) and $b_3 \equiv b_4$  (mod  $2$) follow automatically from \eqref{cong} and \eqref{sum}.  This gives
 \begin{equation}\label{thisgives}
 \mathcal{T}(M) = \frac{1}{8} \underset{\textbf{a}, \textbf{b} }{\left.\sum\right.^{\prime}} \tilde{\psi}\left( \frac{\textbf{a}}{M}, \frac{\textbf{b}}{M}\right)  \mathcal{H}(\textbf{a}, \textbf{b}) 
 \end{equation}
 where $\sum^{\prime}$ indicates the conditions \eqref{cong} -- \eqref{sum},  $\tilde{\psi}$
 is a smooth function such that
 $$
 \tilde{\psi}(\textbf a, \textbf b) = \psi \Big ( \frac{b_1 + b_2}{2} , \frac{b_2 - b_1}{2}, \frac{a_4 - a_3}{2}, \frac{a_1 + a_2}{2}, \frac{a_2 - a_1}{2}, \frac{b_4 - b_3}{2} \Big ) = \psi(\textbf x, \textbf y)
 $$
 and
 \begin{displaymath}
\begin{split}
 \mathcal{H}(\textbf{a},& \textbf{b}) :=     \int_{\Bbb{R}^3}   e\left(\frac{u}{2}(a_1b_2 + a_2 b_1) + \frac{v}{2}(a_3b_4 + a_4b_3) + \frac{z}{4M^2} (a_1+a_2)(b_1+b_2)\right) \\
 &  \widehat{G}\left(ub_1b_2 + va_3a_4 + \frac{z}{4M^2} (b_1+b_2)^2, ua_1a_2 + vb_3b_4 + \frac{z}{4M^2} (a_1+a_2)^2;z, u, v\right)  \,  \dd(u, v, z).  
 \end{split}
 \end{displaymath}
For notational convenience we introduce the following functions in the variables $a_1, \ldots, a_4, b_1, \ldots, b_4$. Put 
\begin{displaymath}
\begin{split}
&P = \max(|a_1 a_2|, |a_3a_4|, |b_1b_2|, |b_3b_4|), \quad \Delta = a_1a_2a_3a_4 - b_1b_2b_3b_4,\\
& \Delta_1 = a_1a_2b_3a_4 + a_1 a_2 a_3 b_4 - a_1 b_2b_3b_4-b_1a_2b_3b_4, \\
& \Delta_2 = a_1b_2a_3a_4 + b_1 a_2 a_3 a_4 - b_1 b_2a_3b_4-b_1b_2b_3a_4.
\end{split}
\end{displaymath}
Note that \eqref{same} and \eqref{sum} imply $P \not= 0$. 

If $\Delta\not= 0$, we can change variables $$u = \frac{a_3a_4 V - b_3b_4 U}{\Delta}, \quad v = \frac{a_1a_2 U - b_1b_2V}{\Delta}$$ to see that $\mathcal{H}(\textbf{a}, \textbf{b})$ equals 
\begin{equation}\label{defH}
\begin{split}
  \frac{1}{| \Delta|}\int_{\Bbb{R}^3}  & \widehat{G}\left(U + \frac{z(b_1+b_2)^2}{4M^2} , V + \frac{z (a_1+a_2)^2}{4M^2}; z, \frac{a_3a_4 V - b_3b_4 U}{\Delta}, \frac{a_1a_2 U - b_1b_2V}{\Delta} \right)  \\
& e\left(\frac{U\Delta_1}{2\Delta} + \frac{V\Delta_2}{2\Delta}  + \frac{z}{4M^2} (a_1+a_2)(b_1+b_2)\right) 
  \,  \dd(U, V, z). 
 \end{split}
 \end{equation}
 By repeated partial integration in $U$ and $V$ we conclude  from \eqref{schwartz} that 
\begin{equation}\label{boundH}
\begin{split}
\mathcal{H}(\textbf{a}, \textbf{b}) & \ll_A \frac{1}{|\Delta|}   \min\left(1, \frac{|\Delta|}{P} \right) \left(\left(1 + \frac{|\Delta_1|}{|\Delta| + P}\right)\left(1 + \frac{|\Delta_2|}{|\Delta| + P}\right)\right)^{-A} \\
&\ll  \frac{1}{|\Delta| + P}\left(\left(1 + \frac{|\Delta_1|}{|\Delta| + P}\right)\left(1 + \frac{|\Delta_2|}{|\Delta| + P}\right)\right)^{-A}
\end{split}
\end{equation}
for any $A \geq 0$ and all $\textbf{a}, \textbf{b} \ll M$. The final bound remains true for $\Delta = 0$, which we show now.  Suppose without loss of generality that $P = |a_1a_2| \not = 0$ (the other cases are similar). We observe that $\Delta = 0$ implies $b_1b_2\Delta_1 = -a_1a_2 \Delta_2$, in particular $|\Delta_2| \leq |\Delta_1|$. 
In this case we change variables simply by
$$ub_1b_2 + va_3a_4 = \frac{b_1b_2}{a_1a_2}U, \quad ua_1a_2 + v b_3b_4  =  U$$
 and recast $\mathcal{H}(\textbf{a}, \textbf{b})$ as
\begin{equation*}
\begin{split}
  \frac{1}{|a_1a_2|}\int_{\Bbb{R}^3} & \widehat{G}\left(\frac{b_1b_2}{a_1a_2}U + \frac{z(b_1+b_2)^2}{4M^2} , U + \frac{z(a_1+a_2)^2}{4M^2} ; z, \frac{U - b_3 b_4 v}{a_1 a_2}, v \right)  \\
& e\left( \frac{U(b_1a_2 + b_2a_1)}{2a_1a_2} + \frac{v\Delta_1}{2 a_1a_2} +  \frac{z}{4M^2} (a_1+a_2)(b_1+b_2)\right) 
  \,  \dd(U, v, z). 
 \end{split}
 \end{equation*}
Now integration by parts with respect to $v$ confirms \eqref{boundH} again. 

In order to evaluate asymptotically \eqref{thisgives}, we would like to apply a kind of Lipschitz principle and replace the sum over $\textbf{a}, \textbf{b}$ by an integral. This is not directly possible, because $\mathcal{H}(\textbf{a}, \textbf{b})$ is quite oscillatory if $\Delta$ is small. 
Using Diophantine techniques, we will establish in the next section the following crucial result, which basically tells us that the contribution of small values of $\Delta$ can be absorbed in an error term. 

\begin{prop}\label{prop4} Let $0 < \delta \leq 1/2$. For any $\varepsilon > 0$ there exists   $\varepsilon' > 0$ so that  $$\underset{\substack{\textbf{a}, \textbf{b} \ll M\\ |\Delta| \leq M^{4-\delta}\\ \Delta_1, \Delta_2 \ll (|\Delta| + P)^{1+\varepsilon'}}}{\left.\sum \right.^{\prime}} \frac{1}{|\Delta| + P} 
 \ll  M^{2- \delta/128 + \varepsilon}. $$
\end{prop} 
  
Let $\phi$ be a smooth function with support 
on $[1/2, \infty]$ that is 1 on $[1, \infty]$, and  write
\begin{equation}\label{defPhi}
\Phi(\textbf{a}, \textbf{b}) :=  \tilde{\psi}\left( \frac{\textbf{a}}{M}, \frac{\textbf{b}}{M}\right)  \phi\left(\frac{|\Delta|}{M^{4-\delta}}\right) .
\end{equation} 
 From \eqref{thisgives}, \eqref{boundH} and Proposition \ref{prop4} we obtain 
\begin{equation}\label{witherror}
\mathcal{T}(M) = \frac{1}{8} \underset{ \textbf{a}, \textbf{b}   }{\left.\sum\right.^{\prime}} \Phi(\textbf{a}, \textbf{b}) \mathcal{H}(\textbf{a}, \textbf{b})  +O(M^{2-\delta/129}). 
\end{equation}
Note that condition \eqref{same} is now automatic since   it is only violated in the case $\Delta = 0$ (which is not in the support of $\Phi$).  On the support of  $\Phi$ we  obtain 
from \eqref{schwartz} and \eqref{defH} that 
\begin{equation}\label{nabla}
\Big \| \nabla \Big ( \Phi(\textbf{a}, \textbf{b})  \mathcal{H}(\textbf{a}, \textbf{b})\Big ) \Big \| \ll \frac{1}{|\Delta|} \cdot  \frac{M^3}{|\Delta|} 
 \ll M^{-5+2\delta}.
\end{equation}
Using \eqref{sum}, we interpret   the main term in \eqref{witherror} as a sum over $a_1, a_2, a_3, b_1, b_2, b_3$ subject only to \eqref{cong}, and in the following we apply the same interpretation for all functions in $\textbf{a}, \textbf{b}$, such as $\Phi$, $\mathcal{H}$, $P$, $\Delta$, $\Delta_1$ and $\Delta_2$. 
 By a standard application of the Euler-Maclaurin summation formula, using \eqref{nabla}, we see that
$$\frac{1}{8} \underset{ \textbf{a}, \textbf{b}   }{\left.\sum\right.^{\prime}} \Phi(\textbf{a}, \textbf{b}) \mathcal{H}(\textbf{a}, \textbf{b})  = \frac{1}{32} \int_{\Bbb{R}^3}\int_{\Bbb{R}^3}  \Phi(\textbf{a}, \textbf{b}) \mathcal{H}(\textbf{a}, \textbf{b})  \,\dd(a_1, a_2, a_3)  \dd(b_1, b_2, b_3)+ O(M^{6-5+2\delta}). $$
   The following continuous analogue of Proposition \ref{prop4} allows us to remove the cut-off function $\phi$ in the definition \eqref{defPhi}.    
 \begin{lemma}\label{lem1} Let $\delta\leq  1/2$. For any $\varepsilon > 0$ there exists $\varepsilon' > 0$ such that  $$\underset{\substack{\textbf{a}, \textbf{b} \ll M\\ |\Delta| \leq M^{4-\delta}\\ \Delta_1, \Delta_2 \ll (|\Delta| + P)^{1+\varepsilon'}}}{\int \cdots \int } \frac{1}{|\Delta| +P} \,  \dd(a_1, a_2, a_3)  \dd(b_1, b_2, b_3) 
 \ll  M^{2- \delta/4}.$$  \end{lemma}
  We will prove this in the final section. Choosing $\delta = 1/2$, we have now arrived at
  \begin{equation}\label{arrived}
  \mathcal{T}(M) =  \frac{1}{32} \int_{\Bbb{R}^3}\int_{\Bbb{R}^3}   \tilde{\psi}\left( \frac{\textbf{a}}{M}, \frac{\textbf{b}}{M}\right) \mathcal{H}(\textbf{a}, \textbf{b})  \,\dd(a_1, a_2, a_3)  \dd(b_1, b_2, b_3) + O\left(M^{2-\frac{1}{258}}\right).
  \end{equation}
 At this point we invert the change of variables $(\textbf{x}, \textbf{y}) \mapsto (\textbf{a}, \textbf{b})$ in \eqref{step3}, undo the integration \eqref{step2} over $\alpha_1, \alpha_3$ and revert the Fourier inversions \eqref{step1}. This shows that
 \begin{equation}\label{im}
 \frac{1}{32}  \int_{\Bbb{R}^3}\int_{\Bbb{R}^3}  \tilde{\psi}\left( \frac{\textbf{a}}{M}, \frac{\textbf{b}}{M}\right)  \mathcal{H}(\textbf{a}, \textbf{b}) \, \dd(a_1, a_2, a_3)  \dd(b_1, b_2, b_3) = \mathcal{I}(M),
 \end{equation}
 where  
 \begin{equation*}
\begin{split}
\mathcal{I}(M) := & \frac{1}{8} \int_{\mathcal{D}} F({\bm \alpha}) \int_{\Bbb{R}^3} \int_{\Bbb{R}^3}   
\psi\left(\frac{\textbf{x}}{M}, \frac{\textbf{y}}{M}\right) W_1\left(\frac{q_{{\bm \alpha}}(x_1, y_1) - q_{{\bm \alpha}}(x_2, y_2)}{D({\bm \alpha}) }\right)\\
&W_2\left(\frac{q_{{\bm \alpha}}(x_1, y_1) - q_{{\bm \alpha}}(x_3, y_3)}{D({\bm \alpha}) }\right)V\left(\frac{q_{{\bm \alpha}}(x_1, y_1)}{M^2 D({\bm \alpha}) }\right) \, \dd\textbf{x} \, \dd\textbf{y} \, \dd_{\text{hyp}}{\bm \alpha}
\end{split}
\end{equation*} 
is the continuous analogue of $\mathcal{T}(M)$. Here we can drop the function $\psi(\textbf{x}/M, \textbf{y}/M)$ because it is redundant. 
 By a change of variables $$x_j \leftarrow \left(x_j + \frac{y_j}{2\alpha_1}\right) \frac{\alpha_1^{1/2}}{(4\alpha_1 \alpha_2 - 1)^{1/4}}, \quad y_j \leftarrow y_j \frac{(4\alpha_1 \alpha_2 - 1)^{1/4}}{(4\alpha_1)^{1/2}}$$
the arguments of $V,W_1, W_2$ become independent of ${\bm \alpha}$, and we obtain
$$\mathcal{I}(M)  =  \text{{\rm vol}}(F) \int_{\Bbb{R}^3} \int_{\Bbb{R}^3}  V\left(\frac{(x_1^2 + y_1^2)\pi}{M^2}\right) W_1( (x_1^2 + y_1^2 - x_2^2 - y_2^2)\pi) W_2((x_1^2 + y_1^2 - x_3^2 - y_3^2)\pi) \,\dd\textbf{x}\, \dd\textbf{y} .$$
Changing to polar coordinates, we get
 $$\mathcal{I}(M)  =  (2\pi)^3 \text{{\rm vol}}(F)  \int_{[0, \infty)^3} V\left(\frac{r_1^2\pi}{M^2}\right) W_1( (r_1^2 - r_2^2)\pi) W_2((r_1^2 - r_3^2)\pi) r_1r_2r_3 \, \dd\textbf{r}.$$
Next changing $r_j \leftarrow r_j^2 \pi$, this simplifies to 
\begin{displaymath}
\begin{split}
 \mathcal{I}(M)  & =    \text{{\rm vol}}(F)  \int_{[0, \infty)^3} V\left(\frac{r_1}{M^2}\right) W_1( r_1 - r_2) W_2(r_1 - r_3)  \,\dd\textbf{r}\\
    & =  \text{{\rm vol}}(F) M^2   \int_0^{\infty} V(r_1) \int_{-\infty}^{Mr_1} \int_{-\infty}^{Mr_1} W_1(r_2) W_2(r_3) \,\dd\textbf{r}
    \end{split}
    \end{displaymath}
and for $M$ sufficiently large this equals  $\text{{\rm vol}}(F) \hat{V}(0) \hat{W}_1(0) \hat{W}_2(0)M^2$. Combining this with \eqref{arrived} and \eqref{im}, we complete the proof of Proposition \ref{prop1}. 
 
 \section{Proof of Theorem \ref{thm:main}: Diophantine analysis}\label{sec:Dioph}

In this section we prove Proposition \ref{prop4}. For notational convenience we introduce the notation $X \preccurlyeq Y$ to mean $X \ll_{\varepsilon} Y M^{\varepsilon}$ for $\varepsilon > 0$, where the meaning of $\varepsilon$ can change from line to line. Next we observe that the conditions \eqref{cong} -- \eqref{sum} as well as  $|\Delta|$, $P$ and  the set $\{|\Delta_1|, |\Delta_2|\}$ are invariant under the following symmetries:
\begin{itemize}
\item interchanging indices 1 and 2;
\item interchanging indices 3 and 4;
\item interchanging indices 1, 2 with 3, 4;
\item interchanging $a$ with $b$.
\end{itemize}
These 4 involutions generate a $2$-subgroup of $S_8$ of order 16. In particular, without loss of generality we can and will assume that 
\begin{equation}\label{assume}
\begin{split}
&\min(|a_3|, |a_4|, |b_3|, |b_4|) \leq \min(|a_1|, |a_2|, |b_1|, |b_2|),\\
& \max(|a_1|, |a_2|, |b_1|, |b_2|) = |a_1|, \\
& |b_2b_3 - a_2a_3| \geq |a_2a_4 - b_2b_4|.
\end{split}
\end{equation}
Now we put all variables into dyadic intervals and suppose that 
 $A _1 \leq |a_1| \leq  2A_1, \ldots, A_4 \leq |a_4| \leq 2 A_4$, $B_1 \leq |b_1| \leq  2B_1, \ldots, B_4 \leq |b_4| \leq 2 B_4$ with $0 \leq A_1, \ldots, B_4 \ll M$. 
 By \eqref{assume} we have
 \begin{equation}\label{assume1}
   \min(A_3, A_4, B_3, B_4) \ll \min(A_1, A_2, B_1, B_2), \quad \max(A_1, A_2, B_1, B_2) = A_1. 
 \end{equation}
 We write $P_0 := \max(A_1A_2, B_1 B_2, A_3A_4, B_3 B_4) \ll M^2.$  In addition we assume that  
   \begin{equation}\label{delta}
  D \leq |\Delta| \leq 2D, \quad     0 \leq D \leq M^{4-\delta}.
     \end{equation}  Let $\mathcal{N}(\textbf{A}, \textbf{B}, D)$ denote the number of 8-tuples $(a_1, \ldots, a_4, b_1, \ldots b_4)$ subject to these size constraints and 
\begin{equation}\label{2}
\Delta_1, \Delta_2 \preccurlyeq D + P_0
\end{equation}
 as well as the conditions \eqref{cong} -- \eqref{sum}. 
 
 Using \eqref{sum} and a standard divisor bound, we have the trivial bound
 \begin{equation}\label{trivial}
 \mathcal{N}(\textbf{A}, \textbf{B}, D) \preccurlyeq P_0^2   \min(A_3, A_4)\min(B_3, B_4)\leq P_0^3. 
\end{equation}

First we dispose of two easy degenerate cases. 
 
 \subsection{Degenerate case I: $a_1 a_2 a_3 a_4 b_1 b_2 b_3 b_4 = 0$}
Let us first assume that one of the variables, say $a_1$, equals 0, but none of the $b$-variables vanishes.  By a divisor argument we can choose the $b$-variables in $\preccurlyeq D$ ways, and then by \eqref{sum} we have $A_3$ (say) choices for $a_2, a_3, a_4$. This gives a total number of
\begin{equation}\label{a}
\preccurlyeq D  A_3 \ll  (D + P_0)  M
\end{equation}
choices.

Next, if   any 3 of the variables $a_1, \ldots, a_4, b_1, \ldots, b_4$   vanish, then by \eqref{sum} it is easy to see that we have at most 
\begin{equation}\label{b}
  \ll  P_0^{1/2}M^2 \ll \ (D + P_0)^{1/2} M^2
\end{equation}   
    choices for the remaining ones. 

Up to symmetry,    the only remaining case is that $a_1$ and exactly one of $b_2, b_3, b_4$ vanishes (recall \eqref{same}). If $b_2 = 0$, then \eqref{2} gives $a_3 a_4 a_2 b_1  \preccurlyeq P_0 + D$ in non-zero variables. There are $ \preccurlyeq P_0 + D $ choices by the usual divisor argument, and after choosing $b_3$, $b_4$ is determined by   \eqref{sum}, which matches the contribution in \eqref{a}. If $b_3 = 0$ (the case $b_4 = 0$ is similar), then \eqref{2} and the first equations in \eqref{sum} give $$P_0 + D \succcurlyeq  a_2a_3 a_4 b_1 - a_3 b_1 b_2 b_4  = a_2a_3b_1(a_4 - b_2).$$
    
If $a_4 \not= b_2$, then we fix $a_2, a_3, b_1, a_4 - b_2$ in $\preccurlyeq D + P_0$ ways as well as $a_4$ in $\ll M$ ways, then the rest is determined by \eqref{sum}, so that we end up with a contribution as in \eqref{a}. 
On the other hand, if $a_4 = b_2$, then everything is determined from $a_2, a_3, a_4$, and we obtain $ P_0^{1/2} M^2 \ll  (D + P_0)^{1/2} M^2$ solutions as in \eqref{b}. 

We summarize that the number $\mathcal{N}_0(\textbf{A}, \textbf{B}, D)$ of tuples where one of the variables vanish is at most
\begin{equation}\label{n0}
\mathcal{N}_0(\textbf{A}, \textbf{B}, D) \preccurlyeq (D + P_0) M + (D + P_0)^{1/2} M^2.
\end{equation}
From now   on we focus on $\mathcal{N}_{\ast}(\textbf{A}, \textbf{B}, D)$ where  $a_1a_2a_3a_4b_1b_2b_3b_4 \not= 0$.

 \subsection{Degenerate case II: $a_2 a_4 = b_2 b_4$}
In this subsection we treat another degenerate case, namely the case $$a_2a_4 = b_2b_4.$$
We call this contribution $\mathcal{N}_{\ast, 0}(\textbf{A}, \textbf{B}, D)$. 
Fixing $b_2, b_4$ (non-zero), the values for $a_2, a_4$ are prescribed up to a divisor function. Moreover, using \eqref{sum} we have $$ \frac{D}{b_2b_4} \gg d := a_1a_3 - b_1b_3= a_1(b_2 - a_4) + a_3(b_4 - a_2) - (a_4-b_2)(a_2 - b_4). $$  
Fix a value for $d$ and suppose initially $d\not= 0$.  Then automatically   $(b_2-a_4, b_4-a_2) \not= (0, 0)$.  Choosing a suitable one of $a_1$ or $a_3$, the other one is determined, and then also $b_1, b_3$ from \eqref{sum}. So the total number of choices is
$$\preccurlyeq B_2B_4  \frac{D}{B_2B_4} (A_1+A_3) \ll D M  \ll (D+P_0) M.$$
 
 Now suppose that $d = 0$, so that $a_1a_3 = b_1b_3$. Put $r = (a_1, b_1)$ and $t = a_1/r$, $s = b_1/r$. Then $(s, t) = 1$ and $a_1 = rt$, $b_1 = rs$. Since $(s, t) = 1$, we must have $s \mid a_3$, say $a_3 = su$. Then $b_3 = tu$. Similarly, we can parametrize the equation $b_2b_4 = a_2a_4$ as $b_2 = xy$, $b_4 = zw$, $a_2 = xz$, $a_4 = yw$ with $(y, z) = 1$. Then \eqref{sum} becomes
 $$rt + xz = tu + zw, \quad su + yw = rs + xy,$$
 in other words
 $$t(r-u) = z(w-x), \quad s(r-u) = y(w-x).$$
We have $ r = u$ if and only if $ w=x$, which is equivalent to the opposite of \eqref{same}, so that this case is excluded.  On the other hand, if $r-u \not= 0 \not= w-x$, then by coprimality $s = t$, $z = y$, so that $$a_1 = b_1, \quad a_2 = b_2 \quad a_3 = b_3, \quad a_4 = b_4.$$ Invoking also \eqref{sum}, we find in this case at most
 $$\ll  A_1A_2A_3  \ll P_0A_3 \ll  (D+P_0)M$$
 solutions. We conclude that
 \begin{equation}\label{nast0}
 \mathcal{N}_{\ast, 0}(\textbf{A}, \textbf{B}, D) \preccurlyeq (D+P_0)M.
 \end{equation} 
From now on we assume $a_2a_4 \not= b_2b_4$.

\subsection{Substituting}\label{subst}  
Solving the first inequality in \eqref{delta} for $b_1$ and substituting into the bound for $\Delta_1$ in \eqref{2}, we obtain
\begin{equation*}
(a_2 a_4 - b_2 b_4)(b_2 b_3 - a_2 a_3) \preccurlyeq \left|\frac{b_2}{a_1}\right| (D + P_0) + D \left|\frac{a_2}{a_1}\right|\ll D + P_0
\end{equation*}
by \eqref{assume1}. We write
\begin{equation}\label{kl}
k = a_2 a_4 - b_2 b_4, \quad l = b_2 b_3 - a_2 a_3.
\end{equation}
By \eqref{assume} we have $|k| \leq |l|$, and by the assumption from the previous subsection we have $kl \not= 0$. We assume $K \leq |k| \leq 2K$, $L \leq |l| \leq 2L$ with 
\begin{equation}\label{LK}
   L \geq K \geq 1, \quad   KL \preccurlyeq D + P_0.
 \end{equation}  
 Note that $a_1, b_1$ are determined from the other variables by \eqref{sum}. We call the corresponding contribution $\mathcal{N}_{\ast}(\textbf{A}, \textbf{B}, D, K, L)$.  
 
 Let $d = (a_2, b_2)$ and $a_2 = d a_2'$, $b_2 = d b_2'$. Then we must have $k = dk'$, $l = d l'$ and $a_4 \equiv \overline{a_2'} k'$ (mod $|b'_2|$), $a_3 \equiv  -\overline{a_2'} l'$ (mod $|b_2'|$).

For the remaining part of the argument,  $\rho$ is a parameter such that 
\begin{equation} \label{rho}
0 < \rho < \frac{\delta}{4}
\end{equation}
with $\delta$ as in Proposition \ref{prop4}, and $\eta$ is a parameter with
$$
0 < \eta < \frac{1}{6}.
$$
For later purposes we take the opportunity to dispose of another somewhat degenerate case at this point. 
Suppose that 
\begin{equation}\label{extra1}
B_4 \leq M^{\eta}, \quad K \asymp A_2A_4\quad  \text{ or } \quad  B_3 \leq M^{\eta}, \quad L \asymp A_2A_3.
\end{equation}
The second set of conditions is analogous, so we focus on the first. 
In this case we fix   $d, a_2', b_2', a_4, b_4, l', a_3$, which then determines $k$ and $b_3$. Thus we see by elementary means that 
\begin{equation}\label{easy}
\begin{split}
\mathcal{N}_{\ast}&(\textbf{A}, \textbf{B}, D, K, L)\\
 & \preccurlyeq \sum_d \sum_{a'_2, b'_2} \sum_{a_4, b_4} \sum_{l'} \sum_{a_3 \equiv  -\overline{a'_2} l' \, (\text{mod }|b'_2|)} 1\ll A_2B_2 A_4B_4 \frac{D+P_0}{K} \left(\frac{A_3}{B_2} + 1\right)\\
& \ll (A_3+B_2)B_4(D+P_0)  \ll (D+P_0)M^{1 + \eta}.
\end{split}
\end{equation}
under the condition \eqref{extra1}. This is sufficient for our purposes, and from now we therefore assume that \eqref{extra1} does not hold.

Let $W$ be a smooth non-negative function with support in $[1/2, 3]$ that is 1 on $[1, 2]$.  We now fix $d, a_2', a_3', k', l', a_3, a_4$, then $b_3, b_4$ are determined, and we obtain
 \begin{equation}\label{beforepoisson}
\begin{split}
\mathcal{N}_{\ast}&(\textbf{A}, \textbf{B}, D, K, L)\\
& \leq \sum_d \sum_{B_2 \leq |b_2| \leq 2B_2} \sum_{(a'_2, b'_2) = 1} W\left(\frac{a'_2}{A_2/d}\right)\sum_{k',l'}W\left(\frac{k'}{K/d}\right) W\left( \frac{l'}{L/d}\right) \\
&  
 \sum_{a_3 \equiv  -\overline{a_2'} l' \, (\text{mod } |b_2'|)} \Phi_{a_2', b_2', -l'; A_3, B_3}(a_3) \sum_{a_4 \equiv  \overline{a_2'} k' \, (\text{mod } |b_2'|)} \Phi_{a_2', b_2', k'; A_4, B_4}(a_4)
\end{split}
\end{equation}
  where
 \begin{equation}\label{Phi} 
  \Phi_{a_2', b_2', n; A, B}(x) = W\left(\frac{|x|}{A}\right) W\left(\frac{|a'_2x - n|}{B |b'_2|}\right).
  \end{equation}
  We have
 \begin{equation}\label{dphi} 
 \partial^{j_1}_x \partial_a^{j_2} \partial_n^{j_3}  \Phi_{a, b_2', n ; A, B}(x)  \ll_{j_1, j_2, j_3} \min\left( A, \frac{BB_2}{A_2}\right)^{-j_1} \left(\frac{B|b_2'|}{A}\right)^{-j_2} ( B|b_2'|)^{-j_3} 
  \end{equation}
  for   $j_1, j_2, j_3 \geq 0$, which we apply with $(n, A, B) = (k', A_3, B_3)$ or $(-l', A_4, B_4).$ 
  
In this analysis we could have interchanged the roles  of $a_2, a_3, a_4$ and $b_2, b_3, b_4$ (leaving $|k|$ and $|l|$ invariant), so that without loss of generality we may and do assume 
\begin{equation}\label{a2b2}
  A_2 \geq B_2. 
\end{equation}

We now distinguish two principal cases which require different treatment of the right hand side of  \eqref{beforepoisson}.  

\subsection{Case 1} Let us assume that
 \begin{equation}\label{case1-assump}
 \min(A_1, \ldots, A_4, B_1, \ldots, B_4) \leq M^{1 - \rho}\quad \text{or} \quad  K L \leq (D + P_0)^{1 - \rho}
 \end{equation}

\subsubsection{Application of Poisson summation}

Here we start with an application of Poisson summation in  \eqref{beforepoisson} with respect to the $a_3, a_4$-sum getting
\begin{equation}\label{case1}
\begin{split}
& \sum_d \sum_{(a'_2, b'_2) = 1} W\left(\frac{a'_2}{A_2/d}\right)\sum_{k',l'}W\left(\frac{k'}{K/d}\right) W\left( \frac{l'}{L/d}\right)\sum_{h_4, h_3\in \Bbb{Z}} e\left(\frac{ \overline{a_2'} k' h_4 - \overline{a_2'} l'  h_3}{|b_2'|}\right) \\
& \frac{1}{|b_2'|^2}  \int_{\Bbb{R}^2} \Phi_{a_2', b_2', k' ; A_4, B_4}(x) \Phi_{a_2', b_2', -l'; A_3, B_3 }(y)e \left( \frac{-xh_4-yh_3}{|b_2'|}\right)\, \dd(x, y). 
 \end{split}
 \end{equation}
By partial integration and \eqref{dphi} we can truncate the $h_j$-sum ($j = 3, 4$) at $$h_j \preccurlyeq \frac{B_2/d}{\min(A_j, B_jB_2/A_2)}  =\frac{1}{d}\left(\frac{B_2}{A_j} + \frac{A_2}{B_j}\right) = : \frac{1}{d} H_j,$$
say, at the cost of a negligible error. For future reference we note that the double integral in \eqref{case1} is trivially bounded by
$$\min\left(A_4, \frac{B_4 B_2}{A_2}\right)\min\left(A_3, \frac{B_3 B_2}{A_2}\right).$$

\subsubsection{The diagonal terms}

The total contribution of the $h_3=h_4 = 0$ term, say $\mathcal{N}_{\ast, \text{diag}}(\textbf{A}, \textbf{B}, D, K, L)$, is 
\begin{equation*}
\begin{split}
& \ll \sum_d \sum_{    (a'_2, b'_2) = 1}\sum_{k',l'} \frac{1}{|b_2'|^2} \int_{\Bbb{R}^2} \Phi_{a_2', b_2', k'; A_4, B_4}(x)  \Phi_{a_2', b_2', l'; A_3, B_3}(y)  \, \dd(x, y)\\
 & \preccurlyeq A_2B_2KL  \min\left(\frac{A_4}{B_2}, \frac{B_4  }{A_2}\right)\min\left(\frac{A_3}{B_2}, \frac{B_3}{A_2}\right)    \ll \min(A_4 B_3, A_3 B_4) KL
   \end{split}
 \end{equation*}
 so that
\begin{equation*}
  \mathcal{N}_{\ast, \text{diag}}(\textbf{A}, \textbf{B}, D, K, L) \ll KL M \min(A_1, \ldots, A_4, B_1, \ldots, B_4)
  \end{equation*}
  by the first inequality in \eqref{assume1}.  By our current assumption \eqref{case1-assump}  we obtain
  \begin{equation}\label{diag1} 
  \mathcal{N}_{\ast, \text{diag}}(\textbf{A}, \textbf{B}, D, K, L) \ll (D+P_0)^{1-\rho} M^2 + (D+P_0)M^{2-\rho}.    \end{equation}

\subsubsection{The off-diagonal terms}

Let us now assume  $(h_3, h_4) \not= (0, 0)$ in \eqref{case1} and call this term $\mathcal{N}_{\ast, \text{off}}(\textbf{A}, \textbf{B}, D, K, L)$. We now apply Poisson summation also to $a_2'$. This gives
\begin{equation}\label{poisson}
\begin{split}
&  \sum_d \sum_{  b_2'}\sum_{k',l'} \frac{1}{|b_2'|^3} W\left(\frac{k'}{K/d}\right)W\left( \frac{l'}{L/d}\right)\sum_{\substack{h_3 \preccurlyeq H_3, h_4 \preccurlyeq H_4\\ (h_3, h_4) \not= 0}} \sum_{r_2}  S(k'h_4 - l'h_3, r_2, |b'_2|)\\
  & \quad\quad \int_{\Bbb{R}^3}W\left(\frac{zd}{A_2}\right) \Phi_{z, b_2', k'; A_4, B_4}(x)  \Phi_{z, b_2', - l';A_3, B_3}(y)e \left( \frac{- yh_3-xh_4 - z r_2}{|b_2'|}\right) \, \dd(x, y, z). 
 \end{split}
  \end{equation}
By partial integration with respect to $z$ and \eqref{dphi} we can truncate the $r_2$-sum at
$$r_2 \preccurlyeq R := \frac{B_2}{A_2}+  \frac{A_4}{B_4}+  \frac{A_3}{B_3} $$
at the cost of a negligible error. Again we distinguish several cases. 

\subsubsection{Off-diagonal case 1} Suppose that 
\begin{equation}\label{off-1}
B_2 B_4 \geq A_2 A_4 M^{-\eta} \quad \text{and} \quad B_2 B_3 \geq A_2 A_3 M^{-\eta}.
\end{equation}
In this case we  apply the Weil bound in \eqref{poisson}.
We first treat the degenerate case where $r_2 = k'h_4-l' h_3 = 0$ in which the Kloosterman sum is large (and in which case $h_3h_4\not= 0$). This contributes at most 
 \begin{displaymath}
\begin{split}
&\preccurlyeq\sum_d  \frac{d}{B_2}      \min\left(\frac{KH_4}{d^2}, \frac{L H_3}{d^2}\right)\frac{A_2}{d} \min\left(A_4, \frac{B_4 B_2}{A_2}\right)\min\left(A_3, \frac{B_3 B_2}{A_2}\right)\\
&\preccurlyeq \frac{1}{B_2}   \left( \left(\frac{B_2}{A_4} +\frac{A_2}{B_4}\right)\left(\frac{B_2}{A_3} +\frac{A_2}{B_3}\right) KL\right)^{1/2} A_2 \min\left(A_4, \frac{B_4 B_2}{A_2}\right)\min\left(A_3, \frac{B_3 B_2}{A_2}\right)\\
& \leq  (KL)^{1/2} \left(B_2(B_3B_4)^{1/2} + (A_2A_3B_2B_4)^{1/2}+ (A_2A_4B_2B_3)^{1/2} + A_2(A_3A_4)^{1/2}\right)\\
& \preccurlyeq  (D + P_0)^{1/2} M^{2}.
\end{split}
\end{displaymath}

The remaining terms contribute by Weil's bound
 \begin{equation}\label{weil}
\begin{split}
& \preccurlyeq   \sum_d \frac{d^{3/2}}{B_2^{3/2}} \frac{KL}{d^2}\left(\frac{H_3}{d} + \frac{H_4}{d} + \frac{H_3H_4}{d^2}\right)(1+R) \frac{A_2}{d}\min\left(A_4, \frac{B_4 B_2}{A_2}\right)\min\left(A_3, \frac{B_3 B_2}{A_2}\right)\\
& \ll \frac{KLA_2}{B_2^{3/2}}  \left(\frac{B_2}{A_4} +\frac{A_2}{B_4} +\frac{B_2}{A_3} +\frac{A_2}{B_3}+ \left(\frac{B_2}{A_4} +\frac{A_2}{B_4}\right)\left(\frac{B_2}{A_3} +\frac{A_2}{B_3}\right) \right)\left(1+\frac{B_2}{A_2}+  \frac{A_4}{B_4}+  \frac{A_3}{B_3}\right)\\
&\quad\quad\quad\quad\quad  \min\left(A_4, \frac{B_4 B_2}{A_2}\right)\min\left(A_3, \frac{B_3 B_2}{A_2}\right)\\
& \ll  KLB_2^{1/2}\left(A_4 + A_3 + A_2 + B_4+B_3 + B_2 + \frac{B_4 A_3}{B_3} + \frac{A_2A_3}{B_3} + \frac{A_4B_3}{B_4} + \frac{A_2A_4}{B_4} \right)\\
& \preccurlyeq  (D+P_0) \left(M^{3/2} + B_2^{1/2}  \left(\frac{B_4 A_3}{B_3} + \frac{A_2A_3}{B_3} + \frac{A_4B_3}{B_4} + \frac{A_2A_4}{B_4} \right)\right).
 \end{split}
\end{equation}

So far we haven't used our assumption \eqref{a2b2}, but we use it now. Combining \eqref{a2b2} with our current assumption \eqref{off-1}, we must have 
$B_4 \geq A_4 M^{-\eta}$ and $B_3 \geq A_3 M^{-\eta}$, so that 
  \eqref{weil} is bounded by 
 $$\preccurlyeq (D + P_0) M^{3/2+\eta}(1 + B_2/A_2) \ll (D + P_0)M^{3/2+\eta},$$
which gives in total
\begin{equation}\label{intotal}
\mathcal{N}_{\ast, \text{off}}(\textbf{A}, \textbf{B}, D, K, L) \preccurlyeq (D + P_0)M^{3/2+\eta} + (D + P_0)^{1/2} M^{2}
\end{equation}
in our present case. 

\subsubsection{Off-diagonal case 2} Let suppose now that
 $$B_2 B_4 \leq A_2 A_4 M^{-\eta}, \quad \text{but} \quad B_2 B_3 \geq A_2 A_3 M^{-\eta}.$$
Then necessarily $ A_2A_4 \asymp K$, otherwise $k = a_2a_4 - b_2b_4$ has no solution. Since we are assuming that \eqref{extra1} does not hold, we have $B_4 \geq M^{\eta}$. We claim now that the present assumptions imply that only the terms $h_4=0$, $h_3 \not= 0$ contribute non-negligibly to \eqref{poisson}. Indeed, applying Poisson summation to \eqref{poisson} with respect to $k'$, we see that up to a negligible error the dual sum has length $\preccurlyeq B_4^{-1} + B_2/K \ll M^{-\eta}$, in other words, 
only the central Poisson term survives.  Since
\begin{equation}\label{kloosterman}
\sum_{k\, (\text{mod }|b'_2|)} S(k' h_4 - l'h_3, r_2, |b_2'|) = 0
\end{equation}
unless $h_4 = 0$, we conclude that this 
  forces $h_4 = 0$, so that $h_3 \not= 0$.  
  
  Having established our claim, we apply  Weil's bound for $S(-h_3l', r_2, |b_2'|)$ and bound $\mathcal{N}_{\ast, \text{off}}(\textbf{A}, \textbf{B}, D, K, L)$ in \eqref{poisson} as above by
  \begin{equation}\label{intotal2}
\begin{split}
 &\preccurlyeq \sum_d \frac{KL/d^2}{B_2^{3/2}/d^{3/2}} \frac{H_3}{d} (1+R)\frac{A_2}{d}\min\left(A_4, \frac{B_4 B_2}{A_2}\right)\min\left(A_3, \frac{B_3 B_2}{A_2}\right) \\
& \preccurlyeq   (D+P_0) \left(M^{3/2} + B_2^{1/2} \frac{B_4 A_3}{B_3}  \right) \\
& \ll  (D + P_0)M^{3/2+\eta} (1 + B_2/A_2) \ll  (D + P_0)M^{3/2+\eta}.
 \end{split}
\end{equation}

\subsubsection{Off-diagonal case 3} Of course, the dual situation $B_2B_3 \leq A_2A_3 M^{-\eta}$, but  $B_2 B_4 \geq A_2A_4 M^{-\eta}$ 
is handled similarly to the previous case.

\subsubsection{Off-diagonal case 4} Finally we treat the case $$B_2 B_4 \leq A_2A_4 M^{-\eta} \quad \text{and} \quad B_2B_3 \leq A_2A_3 M^{-\eta}.$$ Here we can apply Poisson in both $k, l$ in \eqref{poisson}. By the same argument as in \eqref{kloosterman} this forces $h_3 = h_4 = 0$, up to a negligible error,  but this case is excluded in the off-diagonal contribution. 

\subsubsection{The final bound}
Collecting \eqref{intotal} and \eqref{intotal2} we have shown that
\begin{equation*}
\mathcal{N}_{\ast, \text{off}}(\textbf{A}, \textbf{B}, D, K, L) \preccurlyeq (D + P_0)M^{3/2+\eta} + (D +   P_0)^{1/2} M^{2}.
\end{equation*}
Combining this with \eqref{diag1}, we see altogether 
\begin{equation} \label{final2}
\mathcal{N}_{*}(\textbf A, \textbf B, D, K, L) \preccurlyeq (D + P_0)^{1 - \rho} M^2 + (D + P_0) M^{2 - \rho} 
\end{equation}
in the present case \eqref{case1-assump}. 

\subsection{Case 2} We now turn to the second case
\begin{equation}\label{special}
   \min(A_1, \ldots, A_4, B_1, \ldots, B_4) \geq M^{1-\rho}, \quad \text{and} \quad KL \geq (D+P_0)^{1 - \rho}. 
\end{equation}

We could proceed in this case as in Case 1, however we will find that the diagonal term is too large. Therefore we need to introduce an additional constraint into our upper bound for $\mathcal{N}_{*}(\textbf A, \textbf B, D, K, L)$. In this case we will obtain additional constraints for the variable $a_4$.

\subsubsection{Creating an additional constraint on $a_4$}

Solving \eqref{kl} for $b_3, b_4$ and substituting this as well as \eqref{sum} into the definition of $\Delta$,   one finds  by brute force algebraic manipulation 
\begin{equation}\label{eq:first}
 \left(a_4 - \frac{1}{2}\left(b_2 + \frac{k}{a_2}\right)\right)^2 =  H(a_2, b_2, a_3, k, l)  -  \frac{b_2}{a_2l}\Delta 
\end{equation}
where
$$H(a_2, b_2, a_3, k, l) = -\frac{a_3b_2k}{l} + \frac{a_3^2k}{l} + \frac{b^2}{4} + \frac{a_3k}{a_2} - \frac{b_2k}{2a_2} + \frac{k^2}{4a_2^2}.$$
Solving instead \eqref{kl} for $a_3, a_4$, we may   exchange the roles of $a_2, a_3, a_4$ and $b_2, b_3, b_4$, so that again without loss of generality we may assume \eqref{a2b2}. 
Using the coarse bound $\sqrt{x+y} = \sqrt{x} + O(\sqrt{|y|})$, this  implies
\begin{equation}\label{a4simple}
a_4 = G_{\pm}(a_2, b_2, a_3, k, l) + O(E^{1/2}), \quad E :=  D/L, 
\end{equation}
where $G_{\pm}(a_2, b_2, a_3, k, l) = \frac{1}{2}\left(b_2 + \frac{k}{a_2}\right) \pm \sqrt{H(a_2, b_2, a_3, k, l)}$. 
For future reference we note that
\begin{equation}\label{dh}
  \partial_{a_2}^{j}  H(a_2, b_2, a_3, k, l)  \ll_{j } M^2 A_2^{-j}  , \quad   \partial_{a_3}^{j} H(a_2, b_2, a_3, k, l)  \ll_j  M^{2-j}. 
\end{equation}
for $j \geq 0$, where we used that \eqref{a2b2} implies $K/A_2 \ll M$. 

 Notice that because of \eqref{eq:first} and the bound $\Big | \frac{b_2 \Delta}{a_2 \ell} \Big | \ll E = D / L$ that if we have $H(a_2, b_2, a_3, k, l) < - C D / L$ for a sufficiently large absolute constant $C > 0$ then there are no solutions for $a_4$.  We can therefore attach a factor $W((a_4 - G_{\pm}(a_2, b_2, a_3, k, l))/E^{1/2})$ to the $a_4$-sum in \eqref{beforepoisson}. Unfortunately, $G_{\pm}$ may be highly oscillatory in $a_2, a_3$ which is problematic for the application of the Poisson summation formula in a moment. However since we are facing a counting problem where all terms are positive we are free to ``enlarge'' $E$ and $G_{\pm}$. By \eqref{special},  \eqref{LK}, \eqref{delta} and \eqref{rho} we have $$\frac{D}{L} \leq \frac{D}{(D + P_0)^{(1-\rho)/2}} \leq D^{1/2 + \rho/2} \leq M^{(4-\delta)(1+ \rho)/2} \leq M^2.$$ We choose a parameter  $0 < \alpha < 1/6$ and write
\begin{equation}\label{e0}
E_0 := \left(\frac{D+1}{L}\right)^{\alpha}   M^{2(1 - \alpha)} \gg E
\end{equation}
 and obtain 
 $$a_4 = \tilde{G}_{\pm}(a_2, b_2, a_3, k, l) + O(E_0^{1/2}) $$
with  
$$  \tilde{G}_{\pm}(a_2, b_2, a_3, k, l) =  \frac{1}{2}\left(b_2 + \frac{k}{a_2}\right) \pm \sqrt{H(a_2, b_2, a_3, k, l)  + E_0 }.$$  
This gives
\begin{equation}\label{beforepoisson1}
\begin{split}
\mathcal{N}_{\ast}(\textbf{A}, \textbf{B}, D, &K, L) \leq \sum_d \sum_{(a'_2, b'_2) = 1} W\left(\frac{a'_2}{A_2/d}\right)\sum_{k',l'}W\left(\frac{k'}{K/d}\right)W\left( \frac{l'}{L/d}\right) \\
&  
\sum_{a_3 \equiv  -\overline{a_2'} l' \, (\text{mod } |b_2'|)} \Phi_{a_2', b_2', - l'; A_3, B_3}(a_3) \sum_{a_4 \equiv  \overline{a_2'} k' \, (\text{mod } |b_2'|)} \tilde{\Phi}_{d, a_2', b_2', a_3, k', l'}(a_4) 
\end{split}
\end{equation}
  with $\Phi$ as in \eqref{Phi} and 
  $$\tilde{\Phi}_{d, a_2', b_2', a_3, k', l' }(x) = W\left(\frac{x}{A_4}\right) W\left(\frac{a'_2x - k'}{B_4b'_2}\right)\tilde{W}\left(\frac{x - \tilde{G}_{\pm}(a'_2d, b'_2d, a_3, k'd, l'd)}{E_0^{1/2}}\right) .$$
for a suitable smooth weight functions $\tilde{W}$. From \eqref{dh} we obtain
\begin{equation}\label{dg}
\begin{split}
 & \partial_{a_2}^{j}   \tilde{G}_{\pm}(a_2, b_2, a_3, k, l)  \ll_{j }  \frac{M^2}{A_2^j E_0^{1/2}}   + \frac{M^{2j}}{A_2^jE_0^{(2j-1)/2}},  \\  &  \partial_{a_3}^{j} \tilde{G}_{\pm}(a_2, b_2, a_3, k, l)  \ll_{j}  \frac{M^{2-j}}{E_0^{1/2}}+ \frac{M^j}{E_0^{(2j-1)/2}}
  \end{split}
\end{equation}
for $j \geq 1.$ We see now why it is useful to have $E_0$ not too small.

\subsubsection{Application of Poisson summation} 

We apply Poisson summation in $a_3, a_4$ to \eqref{beforepoisson1}. This gives an expression very similar to \eqref{case1} except that $\Phi_{a_2', b_2', k' ; A_4, B_4}(x)$ is replaced with $\tilde{\Phi}_{d, a_2', b_2', y, k' , l'}(x)$. Precisely we get
\begin{equation}\label{case22}
\begin{split}
& \sum_d \sum_{(a'_2, b'_2) = 1} W\left(\frac{a'_2}{A_2/d}\right)\sum_{k',l'}W\left(\frac{k'}{K/d}\right) W\left( \frac{l'}{L/d}\right)\sum_{h_4, h_3\in \Bbb{Z}} e\left(\frac{ \overline{a_2'} k' h_4 - \overline{a_2'} l'  h_3}{|b_2'|}\right) \\
& \frac{1}{|b_2'|^2}  \int_{\Bbb{R}^2} \tilde{\Phi}_{d, a_2', b_2', y, k', l'}(x) \Phi_{a_2', b_2', -l'; A_3, B_3 }(y)e \left( \frac{-xh_4-yh_3}{|b_2'|}\right)\, \dd(x, y). 
 \end{split}
 \end{equation}
In this case, it follows from \eqref{dg} and \eqref{dphi} and integration by parts with respect to $x$ and $y$ in
$$\int_{\Bbb{R}^2} \tilde{\Phi}_{d, a_2', b_2', y, k' , l'}(x) \Phi_{a_2', b_2', -l'; A_3, B_3 }(y)e \left( \frac{-xh_4-yh_3}{|b_2'|}\right)\, \dd(x, y)$$
that we can truncate the $h_j$-sum ($j = 3, 4$) at 
\begin{equation}\label{htilde}
\begin{split}
&h_3 \preccurlyeq \frac{B_2/d}{\min(A_3, B_3B_2/A_2, E_0/M)}  =\frac{1}{d}\left(\frac{B_2}{A_3} + \frac{A_2}{B_3} + \frac{B_2 M}{E_0}\right) = : \frac{1}{d} \tilde{H}_3,\\
&  h_4 \preccurlyeq \frac{B_2/d}{\min(A_4, B_4B_2/A_2, E_0^{1/2})}  =\frac{1}{d}\left(\frac{B_2}{A_4} + \frac{A_2}{B_4} + \frac{B_2}{E_0^{1/2}}\right) = : \frac{1}{d} \tilde{H}_4.
\end{split}
\end{equation}

As before, we denote by $\mathcal{N}_{\ast, \text{diag}}(\textbf{A}, \textbf{B}, D, K, L)$ the contribution of the diagonal terms $(h_3, h_4) = (0,0)$ and by $\mathcal{N}_{\ast, \text{off}}(\textbf A, \textbf B, D, K, L)$ the contribution of the off-diagonal terms $(h_3, h_4) \neq (0,0)$. 

\subsubsection{The diagonal terms}

The introduction of $\tilde{\Phi}$ shortens the $a_4$-sum a bit and by the same argument as in the treatment of diagonal terms in the previous Case 1 together with \eqref{e0} and \eqref{LK}, the central term $h_3=h_4=0$ in \eqref{case22} is bounded by 
\begin{equation}\label{diag2}
\begin{split}
   \mathcal{N}_{\ast, \text{diag}}&(\textbf{A}, \textbf{B}, D, K, L) \ll A_2B_2  \min\left(\frac{\sqrt{E_0}}{B_2}, \frac{A_4}{B_2}, \frac{B_4  }{A_2}\right)\min\left(\frac{A_3}{B_2}, \frac{B_3}{A_2}\right)  KL \\
  & \leq M^{ 1-\alpha} (D+1)^{\alpha/2} B_3  KL^{1- \alpha/2}   \preccurlyeq M^{2-\alpha} (D+1)^{\alpha/2} (D +P_0)^{1 - \alpha/4} \\
  &\leq M^{2-\alpha} (D+1)^{\alpha/4} (D +P_0).\end{split}
\end{equation}

\subsubsection{The off-diagonal terms} 

It remains to obtain an acceptable bound for the off-diagonal contribution of $(h_3, h_4) \neq (0,0)$. 


Because of the condition \eqref{special} all variables have roughly the same size, so this case is much less delicate. As in \eqref{poisson} of Case 1 we apply in \eqref{case22} Poisson summation in the $a_2'$-variable. By \eqref{dg} the dual variable can be truncated at
$$r_2 \preccurlyeq\tilde{R} :=  \frac{B_2}{A_2}+  \frac{A_4}{B_4}+  \frac{A_3}{B_3} + \frac{B_2M^2}{A_2 E_0}.$$
We now run the same computation as already done in Case 1. The terms with $r_2 = kh_4 -lh_3 = 0$ contribute 
 \begin{displaymath}
\begin{split}
\ll &\sum_d  \frac{d}{B_2}      \min\left(\frac{K\tilde{H}_4}{d^2}, \frac{L \tilde{H}_3}{d^2}\right)\frac{A_2}{d} \min\left(A_4, \frac{B_4 B_2}{A_2}\right)\min\left(A_3, \frac{B_3 B_2}{A_2}\right)\\
\preccurlyeq & M^2 (KL)^{1/2} (H_3H_4)^{1/2} \ll M^2 (D + P_0)^{1/2} \left(M^{\rho} + \frac{M^2}{E_0}\right)
 \end{split}
\end{displaymath}
by \eqref{a2b2} and \eqref{special}. Similarly, the remaining terms contribute
\begin{displaymath}
\begin{split}
& \preccurlyeq   \sum_d \frac{d^{3/2}}{B_2^{3/2}} \frac{KL}{d^2}\left(\frac{\tilde{H}_3}{d} + \frac{\tilde{H}_4}{d} + \frac{\tilde{H}_3\tilde{H}_4}{d^2}\right)(1+\tilde{R}) \frac{A_2}{d}\min\left(A_4, \frac{B_4 B_2}{A_2}\right)\min\left(A_3, \frac{B_3 B_2}{A_2}\right)\\
& \preccurlyeq (D + P_0)M^{3/2}(\tilde{H}_3 + \tilde{H}_4 + \tilde{H}_3\tilde{H}_4) (1 + \tilde{R})   \preccurlyeq (D + P_0)M^{3/2}\left(M^{\rho} + \frac{M^2}{E_0}\right)^3. 
  \end{split}
\end{displaymath}
We have shown
$$\mathcal{N}_{\ast, \text{off}}(\textbf{A}, \textbf{B}, D, K, L) \preccurlyeq M^2 (D + P_0)^{1/2} \left(M^{\rho} + \frac{M^2}{E_0}\right) + (D + P_0)M^{3/2}\left(M^{\rho} + \frac{M^2}{E_0}\right)^3. 
$$
 By \eqref{e0} we have 
$$ M^2/E_0 \leq L^{\alpha} M^{2-2(1-\alpha)} \ll M^{4\alpha} $$
and we recall $M \ll P_0^{1/(2-2\rho)} \leq (D + P_0)^{1/(2 - 2\rho)}$. Hence in total we obtain
\begin{equation}\label{offfinal1}
\begin{split}
\mathcal{N}_{\ast, \text{off}}(\textbf{A}, \textbf{B}, D, K, L)& \preccurlyeq (D + P_0)M^{3/2 + \max(3\rho, 12\alpha)}  + (D +   P_0)^{\frac{1}{2} + \frac{\max(\rho, 4\alpha)}{2 - 2\rho}} M^{2}.  
\end{split}
\end{equation}

\subsubsection{The final bound}
Combining \eqref{diag2} and \eqref{offfinal1}, we conclude under the present assumption \eqref{special} that 
\begin{equation} \label{final3}
\begin{split}
\mathcal{N}_{*}(\textbf A, \textbf B, D, K, L) &\preccurlyeq M^{2-\alpha}  (D+1)^{\alpha/4} (D +P_0) \\
&+ (D + P_0)M^{3/2 + \max(3\rho, 12\alpha)}   + (D +   P_0)^{\frac{1}{2} + \frac{\max(\rho, 4\alpha)}{2 - 2\rho}} M^{2}
\end{split}
\end{equation}
for a parameter $0 < \alpha < \frac{1}{6}$ that is at our disposal.

\subsection{The endgame} We summarize our findings. Suppose that
\begin{equation}\label{para}
\max(3\rho, 12\alpha) \leq \frac{1}{2} - \rho, \quad  \frac{\max(\rho, 4\alpha)}{2 - 2\rho} \leq \frac{1}{2} - \rho. 
\end{equation}
Combining \eqref{trivial}, \eqref{n0}, \eqref{nast0}, \eqref{easy}, \eqref{final2} and \eqref{final3}
and summing over $(\log M)^2$ values of $K, L$, we obtain 
\begin{displaymath}
\begin{split}
\mathcal{N}(\textbf{A}, \textbf{B}, D) \preccurlyeq  (D+P_0) \min\left(P_0^2,   M^{2-\rho} + (D + P_0)^{-\rho} M^2 + M^{2-\alpha} (D+1)^{\alpha/4}   \right)
\end{split}
\end{displaymath}
in all cases, where $0 < \rho,   \alpha \leq 1/6$  are subject to \eqref{rho} and \eqref{para}, but otherwise   at our disposal. By \eqref{delta}, we can simplify
\begin{displaymath}
\begin{split}
\mathcal{N}(\textbf{A}, \textbf{B}, D) & \preccurlyeq  (D+P_0) M^2 \min\left((P_0/M)^2,   M^{-\rho} + P_0^{-\rho}  + M^{-\alpha \delta/4}   \right).
\end{split}
\end{displaymath}
Since $\min((P_0/M)^2, P_0^{-\rho}) \leq M^{-2\rho/(2+\rho)}$, we obtain finally
$$\mathcal{N}(\textbf{A}, \textbf{B}, D)  \preccurlyeq  (D+P_0) M^2 \left(M^{-2\rho/(2+\rho)}  + M^{-\alpha \delta/4}   \right).$$
Choosing $\rho = \delta/4 \leq 1/8$, $\alpha = 1/32$, both conditions \eqref{rho} and \eqref{para} are satisfied, and Proposition \ref{prop4} follows. 

 \section{Conclusion of the proof of Theorem \ref{thm:main}: The continuous case} \label{sec:cont}
 
 It is now a simple matter to prove Lemma \ref{lem1}, which follows the strategy of the previous section, but without any arithmetic input. We symmetrize the $6$-dimensional  integration region and consider all 8 variables $a_1, \ldots, a_4$, $b_1, \ldots, b_4$ subject to \eqref{sum}. We put all variables into dyadic intervals $A_j \leq |a_j| \leq 2A_j$, $B_j \leq |b_j| \leq 2B_j$, $D \leq |\Delta| \leq 2D$, $\Delta_1, \Delta_2 \preccurlyeq \Delta  + P_0$ where the meaning of $P_0$ is the same as in the previous section. The numbers $A_j, B_j, D$ run through logarithmically many positive and negative powers of 2 and are bounded by $M^{-100} \leq A_j, B_j  \ll M$, $M^{-100} \leq D \ll M^{4-\delta}$.  (If one of these is $\ll M^{-100}$, the coarsest trivial estimates suffice.)   The same substitution as in Subsection \ref{subst} yields
 $$kl := (a_2a_4 - b_2b_4)(b_2b_3- a_2a_3) \preccurlyeq D + P_0.$$
We also put $k, l$, into logarithmically many dyadic intervals $K \leq |k| \leq 2K$, $L \leq |l| \leq 2L$, $KL \preccurlyeq D+P_0$, and assume $K \leq L$. We also have the relation \eqref{a4simple}, so that the region of integration for $a_4$ is of length $E^{1/2} = (D/L)^{1/2}.$
The total volume of this region can be computed, of course, just as the central Poisson term. We integrate over $a_4, b_4, b_3, a_3, b_2, a_2$ in this order and see that it    is bounded by
\begin{displaymath}
\begin{split}
&\preccurlyeq A_2 B_2     \frac{L}{A_2} B_3  \frac{K}{B_2} \sqrt{E} = KL^{1/2} D^{1/2} B_3 \\
&\ll (D + P_0)^{3/4} D^{1/2} M \leq (D + P_0)  D^{1/4} M \ll (D+P_0) M^{2 - \delta/4}.
\end{split}
\end{displaymath}
Lemma \ref{lem1} follows. 

 \section{Proof of Corollary \ref{cor:liminf}} \label{sec:corliminf}


The main purpose of this section is to deduce Corollary \ref{cor:liminf} from Theorem \ref{thm:main}. The following general lemma implies Corollary \ref{cor:liminf}, but it is stated in such generality that it can also be used in the proof of Corollary \ref{cor:gaps} in Section \ref{sec2}. We prove only the $\liminf$ part of the claim since the proof for the $\limsup$ is identical (up to reversing inequalities). 

\begin{lemma}\label{lem4}
  Let $\sigma$ be a measure with bounded support in $\Bbb{R}$ and let $\gamma$ be a real number. 
  Then for almost all ${\bm \alpha}$ with respect to $\mu_{\text{{\rm hyp}}}$  we have 
  \begin{equation*}
    \liminf_{N \to \infty} \int_{\gamma}^{\infty} T_{3}({\bm \alpha}; [\gamma, x]; N) \dd\sigma(x) \leq \int_{\gamma}^{\infty} (x - \gamma)^2 \dd \sigma(x).
  \end{equation*}
\end{lemma}
\begin{remark}
  The lemma implies a bound for $T_{3}({\bm \alpha}; I; N)$ as in the conclusion of Corollary \ref{cor:liminf}, by writing $I = [a,b]$, picking $\gamma = a$, and setting $\sigma = \delta_{b}$ with $\delta_{b}$ denoting the Dirac measure centered on $b$. Then the right-hand side of Lemma \ref{lem4} becomes $(b-a)^2 = \mu(I)^2$, and thus Corollary \ref{cor:liminf} follows from the lemma.
\end{remark}

\begin{proof}
Let $\varepsilon>0$ be arbitrary. Let $B = B(\varepsilon, \gamma, \sigma)$ denote the set of ${\bm \alpha}$ such that for all $N \geq N_{0}(\varepsilon, \gamma, \sigma)$ we have
\begin{equation}\label{21}
\int_{\gamma}^{\infty} T_3({\bm \alpha}; [\gamma, x]; N)  \dd \sigma(x) > \int_{\gamma}^{\infty}(x - \gamma)^2 \dd \sigma(x) + \varepsilon. 
\end{equation}  We need to show that $\mu_{\text{hyp}}(B) = 0$. Suppose the contrary, i.e.\  that $\mu_{\text{hyp}}(B) > 0$ for some $\varepsilon>0$. For some small $0 < \varepsilon_1 < 1/2$ to be determined later we
choose a (finite) rectangle $R$ with $\mu_{\text{hyp}}(R) > 0$ and $\mu_{\text{hyp}}(B\cap R)  \geq (1-\varepsilon_1)\mu_{\text{hyp}}(R)$. Such an $R$ exists by the Lebesgue density theorem. By Theorem \ref{thm:main} we have
\begin{align*}
\int_R \Big ( \int_{\gamma}^{\infty} T_3({\bm \alpha}; [\gamma, x] ; N) \dd\sigma(x) \Big ) \dd_{\text{hyp}}{\bm \alpha} & = \int_{\gamma}^{\infty} \int_{R} T_{3}({\bm \alpha}; [\gamma, x]; N) \dd_{\text{hyp}}{\bm \alpha} \,\dd \sigma(x) \\ & \leq 
  \int_{\gamma}^{\infty} (x - \gamma)^2 \dd \sigma(x) \mu_{\text{hyp}}(R) + \varepsilon_1
  \end{align*}
  for $N \geq N_1 = N_1(\varepsilon_1, \gamma, \sigma, R)$, so that for all such $N$
   there exists a set $S_{N} \subseteq R$ such that $\mu_{\text{hyp}}(S_N) \geq 2\varepsilon_1\mu_{\text{hyp}}(R)$ and 
\begin{equation}\label{t3}
  \int_{\gamma}^{\infty} T_3({\bm \alpha}; [\gamma, x]; N) \dd\sigma(x) \leq \frac{1}{1 - 2 \varepsilon_1} \Big ( \int_{\gamma}^{\infty} (x - \gamma)^2 \dd \sigma(x) + \varepsilon_1 \Big ) 
  \end{equation}
for all ${\bm \alpha} \in S_N$ and $N \geq N_1$. Let $$S := \limsup S_N = \bigcap_{K \geq 1} \bigcup_{N \geq K} S_N \subseteq R.$$
The set $S$ is the subset of those $\bm \alpha \in R$ which belong to infinitely many $S_{N}$. 
It follows from Lebesgue's dominated convergence theorem that $\mu_{\text{hyp}}(S) \geq 2\varepsilon_1 \mu_{\text{hyp}}(R)$, so that $\mu_{\text{hyp}}(S \cap B) \geq \mu_{\text{hyp}}(S) - \mu_{\text{hyp}}(R \setminus (B \cap R)) \geq \varepsilon_1 \mu_{\text{hyp}}(R)$ and in particular $S \cap B \not= \emptyset$.  On the other hand, if
$$
\int_{\gamma}^{\infty} (x - \gamma)^2 \dd \sigma(x)  + \varepsilon \geq \frac{1}{1 - 2 \varepsilon_1} \Big ( \int_{\gamma}^{x} (x - \gamma)^2 \dd \sigma(x) + \varepsilon_1 \Big )
$$
 then $S \cap B = \emptyset$ by \eqref{21}, and this can be achieved by choosing $\varepsilon_1$ sufficiently small. This contradiction shows $\mu_{\text{hyp}}(B) = 0$. 
\end{proof}


\section{Proof of Corollary \ref{cor:gaps}}\label{sec2}

In this section, we show how Corollary \ref{cor:gaps} can be deduced from Corollary \ref{cor:liminf}. We assume that Theorem \ref{thm:main} holds, and that consequently we may apply Lemma \ref{lem4}. For $b \geq 0$ let
$$  f_{N, {\bm \alpha}}(b) := \frac{1}{N}\big|\{1 \leq i \leq N-1 \mid \Lambda_{i+1}({\bm \alpha}) - \Lambda_i({\bm \alpha})  \leq b\}\big|.$$
This is a piecewise constant non-decreasing function. Let $\varepsilon > 0$ be given. We will show that there is an infinite sequence of indices $N$ for which $f_{N, {\bm \alpha}}(G - \varepsilon) < 1$, with $G$ as in the statement of Corollary \ref{cor:gaps}, thus establishing the existence of infinitely many pairs of consecutive elements of the sequence which are more than $G - \varepsilon$ apart. Suppose to the contrary that $f_{N, {\bm \alpha}}(G - \varepsilon) = 1$ for all large enough $N > N_0({\bm \alpha}, \varepsilon)$. Then,  by partial summation, 
$$\Lambda_{N+1}({\bm \alpha}) - \Lambda_1({\bm \alpha}) =  \sum_{i=1}^{N-1} (\Lambda_{i+1}({\bm \alpha}) - \Lambda_i({\bm \alpha})) =  (N-1) \Big((G - \varepsilon)  - \int_0^{G - \varepsilon} f_{N, {\bm\alpha}} (x) \, \dd x\Big). $$
The left hand side is $(1 + o_{{\bm \alpha}}(1))N$ since the $\Lambda_i(\bm \alpha)$ are normalized so that the average gap in between them is $1 + o(1)$. In other words,
\begin{equation}\label{gn}
G  \geq 1 + (1 - 10^{-6}) \varepsilon + \int_0^{G - \varepsilon} f_{N, {\bm\alpha}} (x) \, \dd x 
\end{equation}
for every  $N \geq N_0 = N_0({\bm \alpha}, \varepsilon)$. 
We now record three elementary inequalities:
\begin{itemize}
\item For all indices $i$ we have
  \begin{equation} \label{first}
  \textbf{1}_{\Lambda_{i + 1} - \Lambda_{i} \leq b}  \geq \sum_{\substack{j > i \\ \Lambda_{j} - \Lambda_{i} \leq b}} 1 - \sum_{\substack{k > j > i \\ \Lambda_{j} - \Lambda_{i} \leq b \\ \Lambda_{k} - \Lambda_{i} \leq b}} 1.
  \end{equation}
\item For all indices $i$ we have
  \begin{equation} \label{second}
  \textbf{1}_{\Lambda_{i + 1} - \Lambda_{i} \leq b} \geq \frac{2}{3} \sum_{\substack{j > i \\ \Lambda_{j} - \Lambda_{i} \leq b}} 1 - \frac{1}{3} \sum_{\substack{k > j > i \\ \Lambda_{j} - \Lambda_{i} \leq b \\ \Lambda_{k} - \Lambda_{i} \leq b }} 1.  
  \end{equation}
\item For all $i > N_0({\bm \alpha}, \varepsilon)$,
  \begin{equation} \label{third}
  \mathbf{1}_{\Lambda_{i + 1} - \Lambda_{i} \leq b} = 1 - \textbf{1}_{\Lambda_{i + 1} - \Lambda_{i} \in (b, G - \varepsilon]} \geq 1 - \sum_{\substack{j > i \\ \Lambda_{j} - \Lambda_{i} \in (b, G - \varepsilon]}} 1.  
  \end{equation}
\end{itemize}
To prove the first two inequalities notice that if there are exactly $\ell \geq 1$ values $\Lambda_{j}$ in the interval $(\Lambda_{i}, \Lambda_{i} + c]$ then the inequalities \eqref{first} and \eqref{second} amount to, respectively,
$$
1 \geq \ell - \binom{\ell}{2} \quad \text{ and }\quad  1 \geq \frac{2}{3} \ell - \frac{1}{3} \binom{\ell}{2}
$$
which are easily seen to hold. 
Since moreover \eqref{first} and \eqref{second} are trivially true when there are no $\Lambda_j$'s in the interval $(\Lambda_{i}, \Lambda_{i} + c]$, we conclude that \eqref{first} and \eqref{second} always hold. The final inequality \eqref{third} follows from our assumption that $\Lambda_{i + 1} - \Lambda_{i} \leq G - \varepsilon$ for all $N > N_0({\bm \alpha}, \varepsilon)$ and from the trivial upper bound
$$
\textbf{1}_{\Lambda_{i + 1} - \Lambda_{i} \in I} \leq \sum_{\substack{j > i \\ \Lambda_{j} - \Lambda_{i} \in I}} 1
$$
that holds for all intervals $I \subset \Bbb{R}$.

Summing \eqref{first}, \eqref{second} and \eqref{third} over $i$ gives, respectively,
\begin{align*}
  f_{N, {\bm \alpha}}(b) & \geq T_{2}(\bm \alpha; [0, b]; N) - \frac{1}{2}  T_{3}({\bm \alpha}; [0,b]; N), \\
  f_{N, {\bm \alpha}}(b) & \geq \frac{2}{3} \cdot T_{2}(\bm \alpha; [0,b]; N) - \frac{1}{6}  T_{3}({\bm \alpha}; [0,b]; N), \\
  f_{N, {\bm \alpha}}(b) & \geq 1 - T_{2}({\bm \alpha}; [b, G - \varepsilon]; N),
\end{align*}
where each of these inequalities holds for all possible values of $b$. Let $C = \sqrt{6 G -5} - 1 \, (\approx 1.65)$. Using the above inequalities we see that
\begin{align} \label{mainineq}
  \int_{0}^{G - \varepsilon} f_{N, {\bm \alpha}}(b) db \geq \int_{0}^{1} & \Big ( T_{2}({\bm \alpha}; [0,b]; N) - \frac{1}{2} T_{3} ({\bm \alpha}; [0,b]; N) \Big ) \dd b \\ \nonumber & + \int_{1}^{C} \Big ( \frac{2}{3}  T_{2}({\bm \alpha}; [0,b]; N) - \frac{1}{6}   T_{3}({\bm \alpha}; [0,b]; N) \Big ) \dd b \\ \nonumber
  & + \int_{C}^{G - \varepsilon} \Big ( 1 - T_{2}({\bm \alpha}; [b, G - \varepsilon]; N) \Big ) \dd b.
\end{align}
We evaluate the integrals using Sarnak's pair correlation result \cite{Sa} (where it is easy to see that the convergence is uniform in intervals $[0, b]$ for $b \leq 1$, say), and using Lemma \ref{lem4} with $\gamma = 0$ and $\sigma$ a measure equal to $\frac{1}{2} \dd x$ on $0 \leq x \leq 1$ resp.\ $\frac{1}{6} \dd x$ on $1 \leq x \leq C$ and vanishing otherwise.  In this way we find that for almost all $\bm \alpha$ there exists a subsequence $N_1 < N_2 < \ldots$ along which the right-hand side of \eqref{mainineq} is
at least
$$
\geq \int_{0}^{1} \Big ( b - \frac{b^2}{2} \Big ) db + \int_{1}^{C} \Big ( \frac{2}{3}   b - \frac{1}{6} b^2 \Big ) db + \int_{C}^{G - \varepsilon} \Big (1 - (G - \varepsilon - b) \Big ) db - 10^{-6} \varepsilon. 
$$
Inserting this lower bound into \eqref{gn} gives
\begin{displaymath}
\begin{split}
G \geq &1 + (1 - 10^{-6}) \varepsilon + \frac{1}{3}   + \Big(3G - \frac{1}{9} \sqrt{6 G - 5} \cdot (3 G + 5)   - \frac{43}{18} \Big)\\
&+  \Big(G   \varepsilon - \frac{G^2}{2}  - \frac{\varepsilon^2}{2} + \sqrt{6 G - 5} \cdot (G - \varepsilon) - 3 G   + 3\Big) - 10^{-6} \varepsilon.
\end{split}
\end{displaymath}
Using the definition of $G$ this simplifies to
$$
0 \geq \big(1 - 2 \cdot 10^{-6} + G - \sqrt{6 G - 5}\big) \varepsilon   - \frac{\varepsilon^2}{2} .
$$
Since $2 \leq G \leq 2.1$, this is a contradiction for all $\varepsilon$ sufficiently small. This proves Corollary \ref{cor:gaps}.

\section{Proof of Theorem \ref{thm4}}\label{sec:proof4}

\subsection{The lower bound}

Let
$$
f_{N}(x) := \frac{1}{N}   \Big \{ i \leq N - 1: \gamma_{i + 1} - \gamma_{i} \leq x \Big \}
\quad \text{and} \quad  
T_{2}(x; N) := \frac{1}{N} \sum_{\substack{\gamma_{i} - \gamma_{j} \in [0, x] \\ 1 \leq j < i \leq N}} 1. 
$$
Suppose without loss of generality that the gaps $\gamma_{i + 1} - \gamma_{i}$ are uniformly bounded by some $G > 0$ for all sufficiently large $i$. Then by integration by parts (as in the previous section), we obtain 
\begin{equation} \label{eq:tocontra}
G \geq 1 + \int_{0}^{G} f_{N}(x) d x + o(1)
\end{equation}
as $N \rightarrow \infty$. 
Notice now that for all $i$ large enough
$$
\mathbf{1}_{\gamma_{i + 1} - \gamma_{i} \leq x} = 1 - \mathbf{1}_{\gamma_{i + 1} - \gamma_{i} \in (x, G]} \geq 1 - \sum_{\substack{\gamma_{i} - \gamma_{j} \in (x, G] \\ 1 \leq j < i \leq N}} 1.
$$
Therefore, summing over $i  < N$ and dividing by $N$ we get
\begin{equation} \label{eq:s}
  f_{N}(x) \geq 1 - (T_{2}(G; N) - T_{2}(x; N)) + o(1)
\end{equation} uniformly in $0 \leq x \leq G$ as $N \rightarrow \infty$.
Given $\varepsilon > 0$ there exists an $N_0(\varepsilon)$ such that for all $N > N_0(\varepsilon)$,
\begin{equation} \label{eq:t}
T_{2}(G;N) - T_{2}(x; N) \leq (G - x) + \varepsilon 
\end{equation}
for all $0 \leq x \leq G$.\footnote{To see this it's enough to split $[0,G]$ into a disjoint union of intervals $I_1, I_2, \ldots$ of length $\varepsilon^2$ and find an $N_0$ large enough so that $T_{2,N}(I_k)$ is less than $\mu(I_k) + \varepsilon^2$ for all $k$.}
In particular fixing $\varepsilon > 0$, combining $f_{N}(x)\geq 0$, \eqref{eq:tocontra}, \eqref{eq:s}, \eqref{eq:t} and taking $N$ to infinity, we conclude that
$$
G \geq 1 + \int_{G - 1}^{G} (1 - G + x) dx + \varepsilon.
$$
Squeezing $\varepsilon$ to zero this implies that $G \geq \frac{3}{2}$. 

\subsection{The upper bound} 

We define a semi-random sequence $(\gamma_i)_{i \geq 1} = (\gamma_i(\omega))_{i \geq 1}$ in the following way. Let $X_1, X_2, \dots$ be a sequence of independent random variables defined on a probability space $(\Omega,\Sigma,\mathbb{P})$, such that $X_i$ has uniform distribution on the interval $[i,i+1]$ for $i \geq 1$. We may assume that none of the $X_i$ has an integer value. Furthermore, for every $m \geq 1$ we define numbers 
$$
y_i^{(m)} = m^2 + \frac{i}{\lceil \sqrt{2m} \rceil}, \qquad 0 \leq i \leq \left\lceil \sqrt{2m} \right\rceil. 
$$
Note that the numbers $y_i^{(m)}, ~0 \leq i \leq \left\lceil \sqrt{2m} \right\rceil$, are all contained in the interval $[m^2,m^2+1]$. Thus these point sets do not intersect or overlap for different values of $m$. Our sequence will consist of the random numbers $(X_i)_i$ from above, as well as of the deterministic numbers $(y_i^{(m)})_{i,m}$. To prevent the random components of the sequence from mixing with the deterministic ones (which would make it unnecessarily complicated to sort the sequence in increasing order), we discard those random numbers whose index is a perfect square, since their values lie in the same range as the values of some of the deterministic numbers. More precisely, we define $(\gamma_i)_{i \geq 1} = (\gamma_i(\omega))_{i \geq 1}$ as the sequence which contains 
\begin{itemize}
 \item all numbers $X_i,~i \geq 1,$ for which $i$ is not a perfect square, as well as 
 \item all numbers $y_i^{(m)}$, for $m \geq 1$ and $0 \leq i < \lceil \sqrt{2 m} \rceil$, 
\end{itemize}
sorted in increasing order. Note that the gaps between consecutive elements of $(\gamma_i)$ are uniformly bounded by $2$.\\

We split the index set $\mathbb{N}$ into two classes $\mathcal{N}_1$ and $\mathcal{N}_2$. The first class $\mathcal{N}_1$ contains all those indices $i$ for which the value of $\gamma_i$ comes from one of the random numbers $(X_j)_{j \geq 1}$. Furthermore, we set $\mathcal{N}_2 = \mathbb{N} \backslash \mathcal{N}_1$, that is, $\mathcal{N}_2$ contains those indices $i$ for which $\gamma_i$ comes from one of the clusters of deterministic points. The set $\mathcal{N}_2$ decomposes into classes $C_m,~m \geq 1$, such that $C_m$ contains those indices $i \in \mathcal{N}_2$ for which $\gamma_i \in [m^2,m^2+1]$. It is easily seen that $|\{i \in \mathcal{N}_2 \mid i \leq N\}| \ll N^{3/4}$. Consequently, 
$$
\left| \left\{ i \in \mathbb{N} \mid \gamma_i \leq N\right\} \right| = N + O (N^{3/4}),
$$
which implies
$$
\frac{\gamma_i}{i} \to 1 \qquad \text{as $i \to \infty$}. 
$$
In other words, the average spacing of $(\gamma_i)$ is 1, as required. We claim that the pair correlation of this sequence is Poissonian, $\mathbb{P}$-almost surely.\\

Assume that $b>0$ is fixed. We split $T_{2}(b; N)$ into different parts in the following way. We write
\begin{eqnarray*}
T_2(b;N) & = & \frac{1}{N} \Big| \Big \{1 \leq i_1 < i_2 \leq N  \mid  \gamma_{i_{2}} - \gamma_{i_{1}} \leq b \Big \}\Big| \\
& =: & \frac{1}{N} \left( {\sum}_1 + {\sum}_2 + {\sum}_3 + \sum_{m=1}^\infty {{\sum}_{4,m}} \right),\\
\end{eqnarray*}
where all sums are taken over indices $1 \leq i_1 < i_2 \leq N$ such that $\gamma_{i_{2}} - \gamma_{i_{1}} \leq b$, subject to the following additional restrictions.
\begin{eqnarray*}
{\sum}_1: & &i_1,i_2 \text{ are both in $\mathcal{N}_1$}. \\
{\sum}_2: & & \text{exactly one of $i_1,i_2$ is in $\mathcal{N}_1$}.\\
{\sum}_3: & & \text{$i_1$ and $i_2$ are both in $\mathcal{N}_2$, but not contained in the same set $C_m$ for some $m$}.\\
{\sum}_{4,m}: & & \text{$i_1$ and $i_2$ are both in $\mathcal{N}_2$ and both contained in the same set $C_m$}.
\end{eqnarray*}
Recall that for $\gamma_{i} \in \mathcal{N}_1$ we have $\gamma_{i} = X_{i} = i + Y_{i}$ where the $Y_{i}$ are independent and uniformly distributed on $[0,1]$. Therefore, 
$$
\frac{1}{N} {\sum}_1 = \sum_{\ell \geq 1} \Big ( \frac{1}{N} \sum_{\substack{k + \ell < N \\ \gamma_{k + \ell}- \gamma_{k} \leq b}} 1 \Big )  = \sum_{\ell \geq 1} \Big ( \frac{1}{N} \sum_{\substack{k + \ell < N \\ Y_{k + \ell} - Y_{k} \leq b - \ell}} 1 \Big ) \to \sum_{\ell \leq b + 1} \mathbb{P}(X - Y \leq b - \ell)
$$
$\mathbb{P}$-almost surely as $N \rightarrow \infty$ by the strong law of large  numbers,   where $X,Y$ are uniformly distributed on $[0,1]$. Since $X,Y$ are independent we have
$$
\mathbb{P}(X - Y \leq x) = \int_{-\infty}^{x} h_{X - Y}(v) \dd v \qquad h_{X - Y}(v) := \int_{\mathbb{R}} \mathbf{1}_{-1 \leq y \leq 0} \cdot \mathbf{1}_{0 \leq v - y < 1} \dd y.
$$
Therefore,
$$
\sum_{\ell \geq 1} \mathbb{P}(X - Y \leq b - \ell) = \int_{-\infty}^{b} \sum_{\ell \geq 1} h_{X - Y}(x - \ell) \dd x.
$$
We now notice that
\begin{align*}
\sum_{\ell \geq 1} h_{X - Y}(x - \ell)  & = \int_{\mathbb{R}} \mathbf{1}_{-1 \leq y \leq 0} \Big ( \sum_{\ell \geq 1} \mathbf{1}_{0 \leq x - \ell - y < 1} \Big ) \dd y \\ & = \int_{\mathbb{R}} \mathbf{1}_{-1 \leq y \leq 0} \cdot \mathbf{1}_{x - y \geq 1} \dd y = \begin{cases} 1 & \text{ if } x \geq 1, \\ x & \text{ if } 0 \leq x \leq 1, \\ 0 & \text{ if } x \leq 0. \end{cases} 
\end{align*}
Combining the above two equations we conclude that 
\begin{equation} \label{sum:1}
 \sum_{\ell \geq 1} \mathbb{P}(X - Y \leq b - \ell) = \begin{cases} 0, & b \leq 0, \\ \frac{b^2}{2}, & b \in [0,1],\\ b - \frac{1}{2}, \quad & b \geq 1. \end{cases} 
\end{equation}
Furthermore, we have
\begin{equation} \label{sum:2}
\frac{1}{N} {\sum}_2 \ll \frac{|\{i \in \mathcal{N}_2 \mid i \leq N\}|}{N} \ll \frac{N^{3/4}}{N} \to 0,
\end{equation}
since for each $i \in \mathcal{N}_2$ there are at most $2 \lceil b \rceil$ indices $j \in \mathcal{N}_1$ such that $|\gamma_{i} - \gamma_{j}| \leq b$. Finally, we have
\begin{equation} \label{sum:3}
\frac{1}{N} {\sum}_3 \to 0 \qquad \text{as $N \to \infty$}, 
\end{equation}
since $b$ is assumed to be fixed and the gaps between $C_m$ and $C_{m+1}$ grow to infinity as $m \to \infty$. Thus it remains to calculate the contribution of the sums ${\sum
}_{4,m}$. So let $m$ be given, and assume that $C_m \subset \{1, \dots, N\}$. Consider 
$$
{\sum}_{4,m} = \Big| \Big \{i_1,i_2 \in C_m, ~i_1 < i_2   \mid  \gamma_{i_{2}} - \gamma_{i_{1}} \leq b \Big \}\Big|.
$$
If $b \geq 1$ then by construction all such pairs of indices are within distance $b$ of each other, and thus for $b \geq 1$ we have
$$
{\sum}_{4,m} = \Big| \Big \{i_1,i_2 \in C_m, ~i_1 < i_2 \Big \}\Big| \sim \frac{|C_m|^2}{2} \sim m.
$$
The maximal $m$ for which $C_m \subset \{1, \dots, N\}$ is of order $\sqrt{N} (1+o(1))$. There may also be one set $C_m$ which is only partially contained in $\{1, \dots, N\}$, but its contribution is of order at most $|C_m|^2 \ll m \ll \sqrt{N}$ and thus negligible. So overall in the case $b \geq 1$ we have
$$
\sum_{m=1}^\infty {{\sum}_{4,m}} \sim \sum_{1 \leq m \leq \sqrt{N}} m \sim \frac{N}{2}.
$$
Thus, together with \eqref{sum:1}, \eqref{sum:2} and \eqref{sum:3}, in the case $b \geq 1$ we obtain
$$
T_2(b;N) \to b - \frac{1}{2} + \frac{1}{2} = b, \qquad \text{as $N \to \infty$, \quad $\mathbb{P}$-almost surely},
$$
which coincides with Poissonian behavior. It remains to consider the case $b \in [0,1]$. So assume that $b$ is in this range and fixed. Let a set $C_m$ be given, and assume again that it is fully contained in $\{1, \dots, N\}$. The points $\gamma_i,~i \in C_m,$ are positioned at the values
$$
m^2, \quad m^2 + \frac{1}{\lceil \sqrt{2m} \rceil}, \qquad m^2 + \frac{2}{\lceil \sqrt{2m} \rceil}, \qquad \dots, \qquad m^2+1.
$$
We could give an explicit formula for the number of pairs of indices $i_1 < i_2$ for which $\gamma_{i_2} - \gamma_{i_1} \leq b$, but it is sufficient to know that the cardinality of this set of pairs is of order 
$$
\sim |C_m|^2 \left((1-b)b + \int_0^b y ~dy \right)  \sim 2m \left(b - b^2/2 \right),
$$
which is easily verified (for example, by looking at the problem in terms of counting lattice points of a square lattice contained in a certain polygon). Thus in the case $b \in [0,1]$ we have
$$
\sum_{m=1}^\infty {{\sum}_{4,m}} \sim \sum_{1 \leq m \leq \sqrt{N}} 2m \left(b - b^2/2 \right) \sim N \left(b - b^2/2 \right).
$$
As a consequence, again in combination with \eqref{sum:1}, \eqref{sum:2} and \eqref{sum:3}, we obtain
$$
T_2(b;N) \to \frac{b^2}{2} + b - \frac{b^2}{2} = b, \qquad \text{as $N \to \infty$, \quad $\mathbb{P}$-almost surely}.
$$
This gives the desired result for $b$ in the range $[0,1]$. Finally we note that by continuity it is sufficient to consider all $b$ contained in a countable dense subset of $[0,\infty)$, such as $\mathbb{Q} \cap [0,\infty)$. For each individual $b$ there is an exceptional set of $\mathbb{P}$-measure zero, and since the union of countably many sets of measure zero also has measure zero, we conclude that there is an element $\omega \in \Omega$ for which $T_2 (b;N) \to b$ holds for all $b \geq 0$. This completes the proof of Theorem \ref{thm4}.

\end{document}